\documentclass{article}
\usepackage{amssymb}
\usepackage{amsmath}
\usepackage{amsthm}

\usepackage{mathrsfs}

\usepackage{cite}

\usepackage{graphicx}
\usepackage{float}

\usepackage{xcolor}

\usepackage{color,hyperref}
\hypersetup{colorlinks,breaklinks,
             linkcolor=red,urlcolor=red,
             anchorcolor=red,citecolor=red}

\usepackage{geometry}
\usepackage[all]{xy}
\input{xy}
\xyoption{all}

\begin{document}

\title{The 2-Adjunction that relates Universal Arrows and Extensive Monads}
\author{\\Adrian Vazquez-Marquez${}^{\dagger}$\footnote{Contact: a.vazquez-marquez@cinvcat.org.mx}
\\Jenylin Zuniga-Apipilhuasco${}^{\ddagger}$
 \\ \\
$\dagger$ Centro de Investigaci\'on en Teor\'ia de Categor\'ias \\
y sus Aplicaciones, A.C., CINVCAT
\\
$\ddagger$ Instituto de Matem\'aticas, UNAM}
\date{}

\maketitle

\begin{center}
\section*{Abstract}
\end{center}

\emph{
In this article the 2-adjunction that relates universal arrows and extensive monads is 
constructed explicitly. This 2-adjunction resembles the one that relates adjunctions and 
monads since the 2-category of universal arrows is isomorphic to the 2-category of adjunctions 
and the 2-category of extensive monads is isomorphic to the 2-category of monads.\\
This article would be useful as a foundation for a theory relating pseudo adjunctions and 
pseudo monads for Gray-categories. On the other hand, it might function as an accesible tool 
for computer scientists on extensive monads.}\\

$\phantom{br}$

\section*{Introduction}\label{2312101714}

It is widely known that adjunctions can be characterized by universal arrows and monads can 
be given as extensive systems or as extension forms \cite{mas-caw-2ned}, \cite{mae-alt}. 
Therefore, the authors take the 2-adjunction that relates adjunctions and monads, as given 
in \cite{clv-kle} or \cite{loj-apk} and rewrite it when the involved 2-categories are 
substituted by their isomorphic counterparts. Hence, the objective of this article is to 
construct, in a very detailed way, this particular 2-adjunction.\\

The authors of this monograph proposed the 2-category of universal arrows as the 2-adjoint 
counterpart for extensive monads based on three \emph{striking} resemblances. For the first 
resemblance, consider a universal arrow $\eta^{RL}$, along this universal arrow there exists 
a function $\zeta^{RL}$ whose image on a morphism of type $v:A\longrightarrow RLB$ is 
$\zeta^{RL}(v):LA\longrightarrow LB$ then $R\zeta^{RL}(v):RLA\longrightarrow RLB$. Furthermore,
if the application for this function is renamed as $(\sim)^{RLB}$ then $v^{RLB}$ 
resembles the behaviour of the respective associated function $(\sim)^{SB}$, for an extensive 
monad, whose image on a morphism $h:A\longrightarrow SB$ is $h^{SB}:SA\longrightarrow SB$.\\

For the second striking resemblace, one of the requirements for an extensive monad
and its associated function is the fulfillment of  the following equation 
$h^{SB}\cdot\eta^{S}\!A = h$. The commutativity of the universal triangle, for the arrow, concludes 
that for $v$, as above, $v^{RLB}\cdot\eta^{RL}\!A = v$.\\

The third and final similarity goes as follows, a classical monad has to have an endofunctor 
$S$ but for an extensive monad this is not the case, $S$  is just a mere function of objects. 
In order to extent this endofunction to morphisms, through the associated extensive function, 
the following formula takes place $Sf := (\eta^{S}B\cdot\! f)^{SB}$. In turn, for the universal 
arrow, if the left function $L$ is to be extended to a left functor, for a morphism $f$, then 
$LF := (\eta^{RL}B\cdot\! f)^{RLB}$ \cite{mas-caw-2ned}.\\

The objetive of this article is two-fold. First, it is to set up the bases for higher-order 
categorical expressions about this relationship. On one hand, F. Marmolejo and R. Wood 
constructed the higher order, and no-iteration, pseudo monad paradigm in their seminal paper 
\cite{maf-moe}. On the other hand, T. M. Fiore proposed a higher order paradigm for 
a universal arrow \cite{fit-psl}. Finally, A. Lauda gives a higher order pseudo adjunction 
that relates adjunctions and pseudo monads. Therefore, this article will complete the 
necessary calculations to complete the whole picture for the higher order pseudo adjuntion.\\

The second objective of the article is to provide a 2-categorical context for computer 
scientists on extensive monads and, its adjoint counterpart, universal arrows.\\

This monograph has two main Theorems, the first one states the 2-adjunction that relates
universal arrows and extensive monads and the second one states the relationship between
this 2-adjunction and the one given by adjunctions and classical monads as in \cite{loj-apk}.\\

The structure of the article is the following. In Section 2, the 2-category of universal arrows is 
constructed. In Section 3, the 2-category of extensive monads is constructed. In Section 4, the 
left 2-functor is defined and in Section 5 its corresponding right 2-functor.\\

In Section 6, the unit and the counit of the pseudo 2-adjunction are defined. In Section 7, the 
2-isomorphisms with adjunctions and classical monads are given.\\

Finally, in Section 8 the conclusions about this monograph are stated. 

\section{Preliminaries and Notation}\label{2310242204}

There is a 2-adjunction between the following 2-categories, the 2-category of adjunctions
where the diagram of right adjoints commutes and the 2-category of Eilenberg-Moore monads
 \cite{clv-kle}, \cite{loj-apk}. In this monograph, the authors prefer to work  with the 
2-category where the right adjoint diagrams no longer commute but they are isomorphic.\\
 
The advantage in taking this approach relies on the generealization of the theory because
one can recover the corresponding to 2-category of adjunctions where the right adjoints 
diagrams commutes. Also, within this approach one can use string diagrams in order to 
represent categorical calculations graphically \cite{joa-get}, \cite{sep-grl},
\cite{tuv-moc}.\\

The notation for the unit of an adjunction, $\mathfrak{L}\dashv\mathfrak{R}$, might be expressed 
as $\eta^{\mathfrak{RL}}$ and the counit as $\varepsilon^{\mathfrak{LR}}$. This notation, in spite of its 
complications, allows one to avoid the proliferation, and the associated confusion, by choosing 
distinct greek letters for units and counits.

\section{The 2-category of Universal Arrows}\label{2311052357}


In this section, the 2-category of universal arrows on the 2-category $Cat$, and denoted 
as $\mathbf{UArr}(Cat)$, is constructed. On the level of objects, it is very well-known that 
universal arrows are equivalent to classical adjunctions \cite{mas-caw-2ned}. This result 
will be extended to any $n$-cell in the 2-category giving an 2-isomorphism with the 2-category 
of classical adjunctions $\mathbf{Adj}_{R}(Cat)$\: \cite{loj-apk}.

\subsection{The Universal Arrows}\label{2311052358}

The objects of this 2-category, the so-called universal arrows, can be given without losing any 
generality as follows. 

\newtheorem{2311060002}{Definition}[subsection]
\begin{2311060002}\label{2311060002}
Given a functor $R:\mathcal{X}\longrightarrow\mathcal{C}$ and a function on objects  
$L:\mathrm{Obj}(\mathcal{C})\longrightarrow \mathrm{Obj}(\mathcal{X})$, a \emph{universal arrow} 
from $A$ to the functor $R$, depicted as $\eta^{RL}\!A:A\longrightarrow RLA$, means that the 
following property holds. For any other morphism $v:A\longrightarrow RX$, in $\mathcal{C}$, 
there exists a unique morphism $w:LA\longrightarrow X$, in $\mathcal{X}$, such that the following 
diagram commutes in $\mathcal{C}$.

\begin{equation}\label{2312180912}
\xy<1cm,1cm>
\POS (16,25) *+{A} = "a12",
\POS (0, 0) *+{RLA} = "a21",
\POS (32, 0) *+{RX} = "a23", 
\POS "a12" \ar_{\eta^{RL}\!A} "a21",
\POS "a12" \ar^{v} "a23",
\POS "a21" \ar_{Rw} "a23"
\endxy
\end{equation}

\begin{flushright}
$\Box$
\end{flushright}

\end{2311060002}

\textbf{Note}: The universal arrow defines a bijection of the following type, for a pair of 
objects $(A,X)$ defined as above, 
$\zeta^{RL}(A,X):\textrm{Hom}_{\mathcal{C}}(A, RX)\longrightarrow \textrm{Hom}_{\mathcal{X}}(LA, X)$
where $\zeta^{RL}(A,X)(v) = w$. This function $\zeta^{RL}$ will be called the 
\emph{associated universal function}, for the corresponding arrow. An alternative short notation 
for this function is the following one $\zeta^{RX}(v) = w$. This notation takes into consideration
only the right-hand vertex at the base.\\

\textbf{Note}: Considera a universal arrow and its universal associated triangle with right-hand
side $v$. For a pair of morphisms $w, w'$ in $\textrm{Hom}_{\mathcal{X}}(LA, X)$ such that their 
$R$-bases make the previous universal triangle commuting then, by unicity, they are equal 
$w = w'$ and so their images under $R$, that is to say $Rw = Rw'$. It can be said that once 
the, parallel, $R$-bases make the same universal triangle commute then these $R$-bases are 
equal.\\

The notation for the complete universal arrow goes as follows 
$(\mathcal{C}, R, \eta^{RL}, \zeta^{R\sim}(\sim))$ or 
$(\mathcal{C}, R, \eta^{RL}, \zeta^{R}(\sim))$.

\subsection{Morphisms of Universal Arrows}\label{2310072324}

The morphisms for universal arrows, denoted as 
$(J,V,\rho^{JV}):(\mathcal{C}, R, \eta^{RL}, \zeta^{R}(\sim))\longrightarrow
(\mathcal{D}, R^{\dagger}, \eta^{RL\dagger}, \zeta^{R\dagger}(\sim))$, 
are comprised of the following structural components.

\begin{enumerate}

\item [i)] A pair of functors $J:\mathcal{C}\longrightarrow\mathcal{D}$ and 
$V:\mathcal{X}\longrightarrow\mathcal{Y}$.

\item [ii)] A natural isomorphism 
$\rho^{JV}:JR\longrightarrow R^{\dagger}V:\mathcal{X}\longrightarrow\mathcal{D}$
depicted by

\begin{equation}\label{2312222042}
\xy
\POS (0,18) *+{\mathcal{C}} = "a11",
\POS (28,18) *+{\mathcal{D}} = "a12",
\POS (0,0) *+{\mathcal{X}} = "a21",
\POS (28,0) *+{\mathcal{Y}} = "a22",
\POS (14,9) *+{\downarrow\rho^{JV}} = "z",
\POS "a21" \ar_{V} "a22",
\POS "a22" \ar_{R^{\dagger}} "a12",
\POS "a21" \ar^{R} "a11",
\POS "a11" \ar^{J} "a12",
\endxy
\end{equation}

\end{enumerate}

\subsection{Transformations for Universal Arrows}\label{2311080007}

The 2-cells for this 2-category are pairs of classical natural transformations. 
For example, 
$(\alpha, \beta):(J, V, \rho^{JV})\longrightarrow (K, W, \rho^{KW}):
(\mathcal{C}, R, \eta^{RL}, \zeta^{R}(\sim))\longrightarrow
(\mathcal{D}, R^{\dagger}, \eta^{RL\dagger}, \zeta^{R\dagger}(\sim))$

\begin{equation*}
\xy<1cm,1cm>
\POS (0,24)  *+{\mathcal{C}} = "a11",
\POS (24,24) *+{\mathcal{D}} = "a12",
\POS (0,0) *+{\mathcal{X}} = "a21",
\POS (24,0) *+{\mathcal{Y}} = "a22",
\POS (10,24) *+{\downarrow\:\alpha} = "x1",
\POS (10,0) *+{\downarrow\:\beta} = "x2",
\POS "a21" \ar^{R}  "a11",
\POS "a22" \ar_{R^{\dagger}} "a12",
\POS "a11" \ar@/^1pc/^{J} "a12",
\POS "a11" \ar@/_1pc/_{K} "a12",
\POS "a21" \ar@/^1pc/^{V} "a22",
\POS "a21" \ar@/_1pc/_{W} "a22",
\endxy
\end{equation*}

\noindent such that 

\begin{equation}\label{2311080011}
R^{\dagger}\beta\circ\rho^{JV} = \rho^{KW}\!\!\circ\alpha R
\end{equation}

This condition can be seen as a surface condition on the previous cylinder.

\subsection{The 2-categorical structure of $\mathbf{UArr}(Cat)$}\label{2501251917}

In order to continue this section a Proposition has to be stated. 

\newtheorem{2501251923}{Proposition}[subsection]
\begin{2501251923}\label{2501251923}
A 2-category $\mathcal{A}$ is a $\textrm{Cat}$-enriched category, therefore for any $A$, $B$ and $C$ 
in $\textrm{Obj}(\mathcal{A})$ there exists a \emph{composition} functor 

\begin{equation}\label{00}
p(A,B,C):\textrm{Hom}_{\mathcal{A}}(A,B)\times\textrm{Hom}_{\mathcal{A}}(B,C)\longrightarrow
\textrm{Hom}_{\mathcal{A}}(A,C)
\end{equation}

The functoriality of this functor is equivalent to the following conditions.

\begin{enumerate}

\item [i)] For any $f$ in $\textrm{Hom}_{\mathcal{A}}(A,B)$, 
$p(A,B,C)[f,\sim]:\textrm{Hom}_{\mathcal{A}}(B,C)\longrightarrow\textrm{Hom}_{\mathcal{A}}(A,C)$
is a functor.

\item [ii)] For any $h$ in $\textrm{Hom}_{\mathcal{A}}(B,C)$, 
$p(A,B,C)[\sim,h]:\textrm{Hom}_{\mathcal{A}}(A,B)\longrightarrow\textrm{Hom}_{\mathcal{A}}(A,C)$
is a functor. 

\item [iii)] For the previous pair 
$p(A,B,C)[f,\sim](h) = p(A,B,C)[\sim,h](f) =: p(A,B,C)[f,h]$.

\item [iv)] For any $\alpha:f\longrightarrow g$ in $\textrm{Hom}_{\mathcal{A}}(A,B)$ and
$\gamma:h\longrightarrow k$ in $\textrm{Hom}_{\mathcal{A}}(B,C)$, the following diagram
commutes. 

\begin{equation}\label{2501252007}
\xy
\POS (0,20) *+{\textrm{Hom}_{\mathcal{A}}(A,C)[f,h]} = "a11",
\POS (60,20) *+{\textrm{Hom}_{\mathcal{A}}(A,C)[f,k]} = "a12",
\POS (0,0) *+{\textrm{Hom}_{\mathcal{A}}(A,C)[g,h]} = "a21",
\POS (60,0) *+{\textrm{Hom}_{\mathcal{A}}(A,C)[g,k]} = "a22",
\POS "a11" \ar_{p(A,B,C)[\alpha,h]} "a21",
\POS "a21" \ar_{p(A,B,C)[g,\gamma]} "a22",
\POS "a11" \ar^{p(A,B,C)[f,\gamma]} "a12",
\POS "a12" \ar^{p(A,B,C)[\alpha,k]} "a22",
\endxy
\end{equation}

\noindent in such a case define 
$p(A,B,C)[\alpha,\gamma]:=p(A,B,C)[g,\gamma]\circ p(A,B,C)[\alpha,h]$.

\end{enumerate}

\end{2501251923}

\textbf{Note}: The morphism $p(A,B,C)[\alpha,h]$ is known as the whiskering of $\alpha$
by $h$.\\

According to the previous Proposition, in order to prove the 2-categorical structure of 
$\mathbf{UArr}(Cat)$ one has to check the functoriality of the  left and right whiskering 
and their respective commutativity \cite{str-cas}.\\

First let us define vertical composition for the 2-category. Consider the following pair of 
2-cells, 
$(\alpha,\beta):(J, V, \rho^{JV})\longrightarrow (F, G, \rho^{FG}):
(\mathcal{C}, R)\longrightarrow(\mathcal{D}, R^{\dagger})$ and
$(\alpha',\beta'):(F, G, \rho^{FG})\longrightarrow(K,W,\rho^{KW}):
(\mathcal{C}, R)\longrightarrow(\mathcal{D}, R^{\dagger})$. The vertical composition 
is defined componentwisely

\begin{equation*}
(\alpha'\circ\alpha,\beta'\circ\beta):(J, V, \rho^{JV})\longrightarrow(K,W,\rho^{KW})
:(\mathcal{C}, R)\longrightarrow(\mathcal{D}, R^{\dagger})
\end{equation*}

It only remains to prove that the condition \eqref{2311080011} is fulfilled.

\begin{equation*}
R^{\dagger}[\beta'\circ\beta]\circ\rho^{JV} 
= R^{\dagger}\beta'\circ R^{\dagger}\beta\circ\rho^{JV} 
= R^{\dagger}\beta'\circ\rho^{FG}\!\!\circ\alpha R
= \rho^{KW}\!\!\circ\alpha' R\circ\alpha R
= \rho^{KW}\!\!\circ[\alpha'\circ\alpha] R
\end{equation*}

$\phantom{br}$

On the other hand, in order to define the whiskering, consider the following
configuration of $n$-cells.

\begin{equation*}
\xy
\POS (0,0) *+{(\mathcal{C},R)} = "a11",
\POS (34,0) *+{(\mathcal{D},R^{\dagger})} = "a12",
\POS (66,0) *+{(\mathcal{Z},R^{\ddagger})} = "a13",
\POS (16,0) *+{\downarrow\:(\alpha,\beta)} = "y",
\POS "a11" \ar@/^1.7pc/^{(J,V,\rho^{JV})} "a12",
\POS "a11" \ar@/_1.7pc/_{(K,W,\rho^{KW})} "a12",
\POS "a12" \ar_{(F,G,\rho^{FG})} "a13",
\endxy
\end{equation*}

The whiskering is defined componentwisely, therefore 
$(F\alpha,G\beta):(FJ,GV,\rho^{FG}V\circ F\rho^{JV})\longrightarrow(FK,GW,\rho^{FG}W\circ F\rho^{KW}):
(\mathcal{C},R)\longrightarrow(\mathcal{Z},R^{\ddagger})$. The condition \eqref{2311080011}
is fulfilled.

\begin{equation*}
R^{\ddagger}G\beta\circ\rho^{FG}V\circ F\rho^{JV} 
= \rho^{FG}W\circ FR^{\dagger}\beta\circ F\rho^{JV}
= \rho^{FG}W\circ F\rho^{KW}\circ F\alpha R
\end{equation*}

The naturality of $\rho^{FG}$ over $\beta$ was used and then the corresponding condition 
over $\alpha$ and $\beta$. The other whiskering is proved in a similar way.\\

Out of these calculations, it is not very hard to check that the horizontal structure
for the following configuration

\begin{equation*}
\xy
\POS (0,0) *+{(\mathcal{C},R)} = "a11",
\POS (34,0) *+{(\mathcal{D},R^{\dagger})} = "a12",
\POS (66,0) *+{(\mathcal{Z},R^{\ddagger})} = "a13",
\POS (16,0) *+{\downarrow\:(\alpha,\beta)} = "y",
\POS (50,0) *+{\downarrow\:(\gamma,\delta)} = "y",
\POS "a11" \ar@/^1.7pc/^{(J,V,\rho^{JV})} "a12",
\POS "a11" \ar@/_1.7pc/_{(K,W,\rho^{KW})} "a12",
\POS "a12" \ar@/^1.7pc/^{(F,G,\rho^{FG})} "a13",
\POS "a12" \ar@/_1.7pc/_{(M,N,\rho^{MN})} "a13",
\endxy
\end{equation*}

\noindent can be well-defined as 

\begin{equation*}
(\gamma,\delta)\ast(\alpha,\beta)
:=(\gamma K\circ F\alpha, \delta W\circ G\beta)
:(FJ,GV,\rho^{FG}V\circ F\rho^{JV})\longrightarrow(MK, NW,\rho^{MN}W\circ M\rho^{KW} )
:(\mathcal{C},R)\longrightarrow(\mathcal{Z},R^{\ddagger})
\end{equation*}

\section{The 2-category of Extensive Monads}\label{2311141919}

In this section, the 2-category of extensive monads on the 2-category $Cat$, denoted as 
$\mathbf{EMnd}(Cat)$, will be given. The respective objects, also known by the name of 
\emph{Kleisli Triples}, were defined by several authors, among them \cite{mae-mos} and
\cite{thw-reb}.\\

It is very well known that extensive monads are equivalent to classical monads 
\cite{mae-mos}. This result will also be extended to any $n$-cell in the 2-category, 
thus giving an 2-isomorphism with the 2-category of classical monads 
$\mathbf{Mnd}_{E}(Cat)$ \cite{loj-apk}.\\

\textbf{Note}: The subscript $\mathbf{Mnd}_{E}(Cat)$ makes reference to the Eilenberg-Moore algebras.

\subsection{Monads in Extensive Form}\label{2311141939}

The objects, or 0-cells, for this 2-category are the \emph{extensive monads}.
An extensive monad, denoted as $(\mathcal{C}, S, \eta^{S},$ $(\sim)^{S})$, consists of the 
following structural components. 

\begin{enumerate}

\item [i)] A function $S:\textrm{Obj}(\mathcal{C})\longrightarrow\textrm{Obj}(\mathcal{C})$.

\item [ii)] A function $\eta^{S}:\textrm{Obj}(\mathcal{C})\longrightarrow\textrm{Mor}(\mathcal{C})$,
such that the image, at $A$, is $\eta^{S}\!A:A\longrightarrow SA$.

\item [iii)] For each pair of objects $A$ and $B$ in $\textrm{Obj}(\mathcal{C})$, a function 
$(\sim)^{SB}:\textrm{Hom}_{\mathcal{C}}(A,SB)\longrightarrow\textrm{Hom}_{\mathcal{C}}(SA,SB)$,
the so-called \emph{associated extensive function} of the monad. The image on 
$h:A\longrightarrow SB$ is denoted as $h^{SB}:A\longrightarrow SB$.

\end{enumerate}

These structural components have to fulfill the following conditions, for any 
$h:A\longrightarrow SB$ and $k:B\longrightarrow SC$.

\begin{enumerate}

\item [a)] $(\eta^{SA})^{SA} = 1_{SA}$,

\item [b)] $h^{SB}\cdot\eta^{SA} = h$,

\item [c)] $(k^{SC}\!\cdot h)^{SC} = k^{SC}\!\cdot h^{SB}$.

\end{enumerate}

\subsection{The algebras for an extensive monad}\label{2311141948}

In order to define morphisms of extensive monads, the category of extensive algebras 
has to be defined first.\\

The objects or algebras for the monad, $(\mathcal{C}, (\sim)^{S})$, are pairs of the 
form $(N, (\sim)^{N})$ where $N$ is in $\mathcal{C}$ and 
$(\sim)^{N}:\textrm{Hom}_{\mathcal{C}}(A, N)\longrightarrow \textrm{Hom}_{\mathcal{C}}(SA, N)$ 
is a function. The following conditions have to be fulfilled for any $a:A\longrightarrow N$, 
$b:B\longrightarrow N$ and $h:A\longrightarrow SB$.\\

\begin{enumerate}

\item [i)] $a^{N}\!\!\cdot \eta^{S}\!A = a$,

\item [ii)] $(b^{N}\cdot h)^{N} = b^{N}\cdot h^{SB}$

\end{enumerate}

\textbf{Note}: These objects might also be called $S$-algebras. 
The object $N$ will be called the algebra and any morphism  $a:A\longrightarrow N$, for example, 
will be called an \emph{algebraic morphism} over, the algebra, $N$.
This can be rephrased as follows, any morphism of type $a:A\longrightarrow N$ is an 
$S$-algebraic morphism.\\

The morphisms of algebras are morphisms $q:(M, (\sim)^{M})\longrightarrow (N, (\sim)^{N})$, 
such that for any $a:A\longrightarrow M$

\begin{equation}
(q\cdot a)^{N} = q\cdot a^{M} 
\end{equation}

This arrange itself into a categorical structure and
the notation for this category is $\mathcal{C}^{S}$.\\

Once the category of algebras for the extensive monad has been constructed, there exists a 
\emph{forgetful} functor $U^{S}:\mathcal{C}^{S}\longrightarrow\mathcal{C}$. On objects it is 
defined as $U^{S}(N,(\sim)^{N}) = N$ and on morphisms as $U^{S}(q) = q$. Also, there exists a 
function on objects $F^{S}\!A = (SA,(\sim)^{SA})$. This is an algebra, since it inherits the 
respective properties for the extensive monad.   

\subsection{Morphisms of Extensive Monads}\label{2311080023}

Resuming with the construction of the 2-category of extensive monads, the morphisms are pairs 
$(P,(\sim)^{PS}):(\mathcal{C},(\sim)^{S})\longrightarrow(\mathcal{D},(\sim)^{T})$ 
with the following structural components.

\begin{enumerate}

\item [i)] $P:\mathcal{C}\longrightarrow\mathcal{D}$ is a functor.

\item [ii)] A function 
$(\sim)^{PS}:\textrm{Hom}_{\mathcal{D}}(D,PSA)\longrightarrow\textrm{Hom}_{\mathcal{D}}(TD,PSA)$,
whose image on $p:D\longrightarrow PSA$ is denoted as $p^{PSA}:TD\longrightarrow PSA$.

\end{enumerate}

These structural components have to fulfill the following conditions.

\begin{enumerate}

\item [i)] $(PSA, (\sim)^{PSA})$ is a $T$-algebra for any $A$ in $\textrm{Obj}(\mathcal{C})$.

\item [ii)] $Ph^{SB}:PSA\longrightarrow PSB$ is a $T$-algebra morphism for any 
$h:A\longrightarrow SB$. Therefore, for $p:D\longrightarrow PSA$ the following condition takes place. 

\begin{equation}
(Ph^{SB}\cdot p)^{PSB} = Ph^{SB}\cdot p^{PSA} 
\end{equation}

\end{enumerate}

\subsection{Transformations of Extensive Monads}\label{2311131247}

The 2-cells in this 2-category are classical natural transformations of type 
$\theta:P\longrightarrow Q:\mathcal{C}\longrightarrow\mathcal{D}$ and denoted as 
$\theta:(P, (\sim)^{PS})\longrightarrow (Q, (\sim)^{QS})$, where the component 
$\theta SA: PSA\longrightarrow QSA$ is required to be a morphism of $T$-algebras, 
that is to say it is a morphism of type
$\theta SA: (PSA,(\sim)^{PSA})\longrightarrow (QSA,(\sim)^{QSA})$.\\

\subsection{The 2-categorical structure of $\mathbf{EMnd}(Cat)$}\label{2312241729}

The proof for the 2-categorical structure will also follow the requirements given 
by Proposition . This proof will also use a pair of Lemmas.\\

Consider the following configuration of 2-cells in $\mathbf{EMnd}(Cat)$.

\begin{equation*}
\xy
\POS (0,0) *+{(\mathcal{C},S)} = "a11",
\POS (36,0) *+{(\mathcal{D},T)} = "a12",
\POS (18,4) *+{\downarrow\theta} = "y",
\POS (18,-7) *+{\downarrow\vartheta} = "y",
\POS "a11" \ar@/^2.6pc/^{(P,(\sim)^{PS})} "a12",
\POS "a11" \ar_{(X,(\sim)^{XS})} "a12",
\POS "a11" \ar@/_2.6pc/_{(Q,(\sim)^{QS})} "a12",
\endxy
\end{equation*}

Since the underlying structure for $\mathbf{EMnd}(Cat)$ is the structure on $Cat$, 
the vertical composition of 2-cells is the usual composition of natural transformations 
in $Cat$.

\newtheorem{2312251949}{Lemma}[subsection]
\begin{2312251949}\label{2312251949}
The vertical compostion is well-defined.
\end{2312251949}

\begin{flushleft}
\emph{Proof}:
\end{flushleft}

Consider the component at $SA$, in $\mathcal{C}$, for the composition 

\begin{equation*}
(\vartheta\circ\theta)(SA) = \vartheta SA\cdot \theta SA.
\end{equation*}

Therefore the composition of components, at $SA$, is a composition of morphisms of 
$T$-algebras which gives a morphism of $T$-algebras. 

\begin{flushright}
$\Box$
\end{flushright}

Consider the following pair of morphisms in $\mathbf{EMnd}(Cat)$, 
$(P,(\sim)^{PS}):(\mathcal{C},S)\longrightarrow(\mathcal{D}, T)$ and
$(W,(\sim)^{WT}):(\mathcal{D},T)\longrightarrow(\mathcal{X}, U)$. Its composition is defined by
the composition of functors $WP$ and the following associated extensive function 
$(\sim)^{W\!PS}:\mathrm{Hom}_{\mathcal{C}}(X,WPSA)\longrightarrow\mathrm{Hom}_{\mathcal{X}}(UX, WPSA)$. 
For $w:X\longrightarrow WPSA$

\begin{equation}\label{2312241907}
w^{W\!PSA} := W(1_{PSA})^{PSA}\cdot[W\eta^{T}PSA\cdot w]^{WTPSA}
\end{equation}

\newtheorem{2312241647}[2312251949]{Lemma}
\begin{2312241647}\label{2312241647}
The horizontal compositon between 1-cells is well-defined.
\end{2312241647}

\begin{flushleft}
\emph{Proof}:
\end{flushleft}

Consider $(P,(\sim)^{PS}):(\mathcal{C},S)\longrightarrow(\mathcal{D}, T)$ and
$(W,(\sim)^{WT}):(\mathcal{D},T)\longrightarrow(\mathcal{X}, U)$ as above, then it has
to be proved that for $A$, in $\mathcal{C}$, $(WPSA,(\sim)^{W\!PSA})$ is a $U$-algebra.
In such a case, for any 
$x:Z\longrightarrow UX$, the following equality has to be fulfilled 
$(w^{W\!PSA}\cdot x)^{W\!PSA} = w^{W\!PSA}\cdot x^{U\!X}$.

\begin{eqnarray*}
(w^{W\!PSA}\cdot x)^{W\!PSA} 
&=& W1_{PSA}^{\:\:PSA}\cdot[W\eta^{T}PSA\cdot W1_{PSA}^{\:\:PSA}\cdot(W\eta^{T}PSA\cdot w)^{WTPSA}\cdot x]^{WTPSA}\\
&=& W1_{PSA}^{\:\:PSA}\cdot[W(\eta^{T}PSA\cdot 1_{PSA}^{\:\:PSA})^{TPSA}\cdot W\eta^{T}PSA\cdot(W\eta^{T}PSA\cdot w)^{WTPSA}\cdot x]^{WTPSA}\\
&=& W1_{PSA}^{\:\:PSA}\cdot W(\eta^{T}PSA\cdot 1_{PSA}^{\:\:PSA})^{TPSA}\cdot[W\eta^{T}PSA\cdot(W\eta^{T}PSA\cdot w)^{WTPSA}\cdot x]^{WTPSA}\\
&=& W(1_{PSA}^{\:\:PSA}\cdot\eta^{T}PSA\cdot 1_{PSA}^{\:\:PSA})^{PSA}\cdot[W\eta^{T}PSA\cdot(W\eta^{T}PSA\cdot w)^{WTPSA}\cdot x]^{WTPSA}\\
&=& W(1_{PSA}^{\:\:PSA})^{PSA}\cdot[W\eta^{T}PSA\cdot(W\eta^{T}PSA\cdot w)^{WTPSA}\cdot x]^{WTPSA}\\
&=& W1_{PSA}^{\:\:PSA}\cdot W1_{TPSA}^{\:\:TPSA}\cdot[W\eta^{T}PSA\cdot(W\eta^{T}PSA\cdot w)^{WTPSA}\cdot x]^{WTPSA}\\
&=& W1_{PSA}^{\:\:PSA}\cdot [W1_{TPSA}^{\:\:TPSA}\cdot W\eta^{T}PSA\cdot(W\eta^{T}PSA\cdot w)^{WTPSA}\cdot x]^{WTPSA}\\
&=& W1_{PSA}^{\:\:PSA}\cdot [(W\eta^{T}PSA\cdot w)^{WTPSA}\cdot x]^{WTPSA}\\
&=& W1_{PSA}^{\:\:PSA}\cdot (W\eta^{T}PSA\cdot w)^{WTPSA}\cdot x^{UX}\\
&=& w^{W\!PSA}\cdot x^{UX}
\end{eqnarray*}

The second equality takes place because of the following property 
$\eta^{T}PSA\cdot 1_{PSA}^{\:\:PSA} = (\eta^{T}PSA\cdot 1_{PSA}^{\:\:PSA})^{TPSA}\cdot \eta^{T}PSA$.
The third equality follows since $W(\eta^{T}PSA\cdot 1_{PSA}^{\:\:PSA})^{TPSA}$ is a morphism of $U$-algebras.
In the fourth equality, that $1_{PSA}$ is an algebraic morphism, over $PSA$, was used.
In the fifth equality, the property $1_{PSA}^{\:\:PSA}\cdot\eta^{T}PSA = 1_{PSA}$ was applied.
For the sixth equality, that $1_{PSA}$ is an algebraic morphism was used again over the morphism $1_{TPSA}$.
The seventh equality follows since $W1_{TPSA}^{\:\:TPSA}$ is morphism of $U$-algebras. 
In the eight equality, the property $1_{TPSA}^{\:\:TPSA}\cdot\eta^{T}PSA = 1_{TPSA}$ was applied.
Finally, in the nineth equality, that $(WTPSA, (\sim)^{WTPSA})$ is a $U$-algebra was used.\\

Now, it has to be proved, for $h:A\longrightarrow SB$, that 
$WPh^{SB}:(WPSA,(\sim)^{W\!PSA})\longrightarrow(WPSB,(\sim)^{W\!PSB})$ is
a morphism of $U$-algebras.

\begin{eqnarray*}
WPh^{SB}\cdot w^{WPSA} 
&=& WPh^{SB}\cdot W1_{PSA}^{\:\:PSA}\cdot[W\eta^{T}PSA\cdot w]^{WTPSA}\\
&=& W(Ph^{SB}\cdot 1_{PSA})^{PSB}\cdot[W\eta^{T}PSA\cdot w]^{WTPSA}\\ 
&=& W(1_{PSB}\cdot Ph^{SB})^{PSB}\cdot[W\eta^{T}PSA\cdot w]^{WTPSA}\\
&=& W(1_{PSB}^{\:\:PSB}\cdot\eta^{T}PSB\cdot Ph^{SB})^{PSB}\cdot[W\eta^{T}PSA\cdot w]^{WTPSA}\\ 
&=& W1_{PSB}^{\:\:PSB}\cdot W(\eta^{T}PSB\cdot Ph^{SB})^{TPSB}\cdot[W\eta^{T}PSA\cdot w]^{WTPSA}\\
&=& W1_{PSB}^{\:\:PSB}\cdot[W(\eta^{T}PSB\cdot Ph^{SB})^{TPSB}\cdot W\eta^{T}PSA\cdot w]^{WTPSB}\\
&=& W1_{PSB}^{\:\:PSB}\cdot[W\eta^{T}PSB\cdot WPh^{SB}\cdot w]^{WTPSB}\\ 
&=& (WPh^{SB}\cdot w)^{WPSB}
\end{eqnarray*}

In the second equality, the property for $Ph^{SB}$ to be a morphism of $T$-algebras was applied.
In the fourth equality, the property $1_{PSB}^{\:\:PSB}\cdot\eta^{T}PSB = 1_{PSB}$ was used. 
In the fifth equality, it was used the property for $1_{PSB}$ being an algebraic morphism, over $PSB$. 
In the sixth equality, the requirement for $W(\eta^{T}PSB\cdot Ph^{SB})^{TPSB}$ being a morphism 
of $U$-algebras was applied. 
Finally, in the seventh equality, the property 
$(\eta^{T}PSB\cdot Ph^{SB})^{TPSB}\cdot \eta^{T}PSA = \eta^{T}PSB\cdot Ph^{SB}$ 
was used.

\begin{flushright}
$\Box$
\end{flushright}

\newtheorem{2312102012}[2312251949]{Proposition}
\begin{2312102012}\label{2312102012}
$\mathbf{EMnd}(Cat)$ has a 2-categorical structure.
\end{2312102012}

\begin{flushleft}
\emph{Proof}:
\end{flushleft}

Consider the following configuration of $n$-cells in $\mathbf{EMnd}(Cat)$

\begin{equation*}
\xy
\POS (0,0) *+{(\mathcal{C},S)} = "a11",
\POS (28,0) *+{(\mathcal{D},T)} = "a12",
\POS (58,0) *+{(\mathcal{X},U)} = "a13",
\POS (12,0) *+{\downarrow\:\theta} = "y",
\POS "a11" \ar@/^1.2pc/^{(P,(\sim)^{PS})} "a12",
\POS "a11" \ar@/_1.2pc/_{(Q,(\sim)^{QS})} "a12",
\POS "a12" \ar_{(W,(\sim)^{WT})} "a13",
\endxy
\end{equation*}

The pull forward composition, or whiskering, is defined as
$W\theta: (WP,(\sim)^{WPS})\longrightarrow(WQ,(\sim)^{WQS}):(\mathcal{C},S)\longrightarrow(\mathcal{X},U)$,
that is to say, with the underlying natural transformation $W\theta$. It is well-defined since 
the component, at $SA$, is a morphism of $U$-algebras, 
$W\theta SA\cdot w^{WPSA} = (W\theta SA\cdot w)^{WQSA}$. 

\begin{eqnarray*}
W\theta SA\cdot w^{WPSA} 
&=& W\theta SA \cdot W1_{PSA}^{\:\:PSA}\cdot[W\eta^{T}PSA\cdot w]^{WTPSA}\\
&=& W(\theta SA\cdot 1_{PSA})^{QSA}\cdot[W\eta^{T}PSA\cdot w]^{WTPSA}\\
&=& W(1_{QSA}^{\:\:QSA}\cdot\eta^{T}QSA\cdot \theta SA)^{QSA}\cdot[W\eta^{T}PSA\cdot w]^{WTPSA}\\
&=& W1_{QSA}^{\:\:QSA}\cdot W(\eta^{T}QSA\cdot \theta SA)^{TQSA}\cdot[W\eta^{T}PSA\cdot w]^{WTPSA}\\
&=& W1_{QSA}^{\:\:QSA}\cdot [W(\eta^{T}QSA\cdot \theta SA)^{TQSA}\cdot W\eta^{T}PSA\cdot w]^{WTQSA}\\
&=& W1_{QSA}^{\:\:QSA}\cdot [W\eta^{T}QSA\cdot W\theta SA\cdot w]^{WTQSA}\\
&=& (W\theta SA\cdot w)^{WQSA}
\end{eqnarray*}

In the second equality, that the component $\theta SA$ is a morphism of $T$-algebras was used.
In the third equality, the property $1_{QSA}^{\:\:QSA}\cdot\eta^{T}QSA = 1_{QSA}$ was applied.
In the fourth equality, that $1_{QSA}$ is an algebraic morphism over the algebra $QSA$ was used.
In the fifth equality, $W(\eta^{T}QSA\cdot \theta SA)^{TQSA}$ being a morphism of $U$-algebras was applied.
Finally, in the sixth equality, the property 
$(\eta^{T}QSA\cdot \theta SA)^{TQSA}\cdot \eta^{T}PSA = \eta^{T}QSA\cdot \theta SA$ was used.\\

Consider the following configuration of $n$-cells in $\mathbf{EMnd}(Cat)$

\begin{equation*}
\xy
\POS (0,0) *+{(\mathcal{C},S)} = "a11",
\POS (30,0) *+{(\mathcal{D},T)} = "a12",
\POS (58,0) *+{(\mathcal{X},U)} = "a13",
\POS (43,0) *+{\downarrow\:\xi} = "y",
\POS "a12" \ar@/^1.2pc/^{(W,(\sim)^{WT})} "a13",
\POS "a12" \ar@/_1.2pc/_{(K,(\sim)^{KT})} "a13",
\POS "a11" \ar_{(Q,(\sim)^{QS})} "a12",
\endxy
\end{equation*}

The pull back composition, or whiskering, is defined as 
$\xi Q:(WQ,(\sim)^{WQS})\longrightarrow(KQ,(\sim)^{KQS}): (\mathcal{C},S)\longrightarrow(\mathcal{X},U)$.
It has to be proved that the component $\xi QSA$ is a morphism of $U$-algebras. Therefore,

\begin{eqnarray*}
\xi QSA\cdot w^{WQSA} 
&=& \xi QSA\cdot W1_{QSA}^{\:\:QSA}\cdot [W\eta^{T}QSA\cdot w]^{WTQSA}\\
&=& K1_{QSA}^{\:\:QSA}\cdot \xi TQSA\cdot [W\eta^{T}QSA\cdot w]^{WTQSA}\\
&=& K1_{QSA}^{\:\:QSA}\cdot [\xi TQSA\cdot W\eta^{T}QSA\cdot w]^{KTQSA}\\
&=& K1_{QSA}^{\:\:QSA}\cdot [K\eta^{T}QSA\cdot \xi QSA\cdot w]^{KTQSA}\\
&=& (\xi QSA\cdot w)^{KQSA}
\end{eqnarray*}

In the second equality, the naturality of $\xi$ over the morphism $1_{QSA}^{\:\:QSA}$ was used.
In the third equality, the fact that $\xi TQSA$ is a morphism of $U$-algebras was applied. 
In the fourth equality, the naturality of $\xi$ over the morphism $\eta^{T}QSA$ was used.\\ 

Define the horizontal composition for the following configuration of 2-cells

\begin{equation}\label{2312261248}
\xy
\POS (0,0) *+{(\mathcal{C},S)} = "a11",
\POS (28,0) *+{(\mathcal{D},T)} = "a12",
\POS (56,0) *+{(\mathcal{X},U)} = "a13",
\POS (13,0) *+{\downarrow\:\theta} = "y",
\POS (41,0) *+{\downarrow\:\xi} = "z",
\POS "a11" \ar@/^1.2pc/^{(P,(\sim)^{PS})} "a12",
\POS "a11" \ar@/_1.2pc/_{(Q,(\sim)^{QS})} "a12",
\POS "a12" \ar@/^1.2pc/^{(W,(\sim)^{WT})} "a13",
\POS "a12" \ar@/_1.2pc/_{(K,(\sim)^{KT})} "a13",
\endxy
\end{equation}

as follows

\begin{equation}\label{2312261249}
\xi\ast\theta := \xi Q\circ W\theta:(WP, (\sim)^{WPS})\longrightarrow(WQ,(\sim)^{WQS})
:(\mathcal{C},S)\longrightarrow(\mathcal{X},U)
\end{equation}

$\phantom{br}$

This is clearly equal to $K\theta\circ\xi P$ as 2-cells in $\mathbf{EMnd}(Cat)$.

\begin{flushright}
$\Box$
\end{flushright}

\section{The Left 2-Functor}\label{2311060013}

In this section, the 2-functor 
$\Phi:\mathbf{UArr}(Cat)\longrightarrow \mathbf{EMnd}(Cat)$ is described. This description 
is developed through each $n$-cell.

\subsection{On Universal Arrows}\label{2311060014}

The 2-functor acts on 0-cells, namely universal arrows, as follows. Consider
$(\mathcal{C},R,\eta^{RL},\zeta^{R}(\sim))$, 
first the induced endofunction for the extensive monad is 
$RL:\textrm{Obj}(\mathcal{C})\longrightarrow\textrm{Obj}(\mathcal{C})$.\\

Second, in order to define the function 
$(\sim)^{RLB}:\textrm{Hom}_{\mathcal{C}}(A, RLB)\longrightarrow \textrm{Hom}_{\mathcal{C}}(RLA, RLB)$,
for the extensive monad $RL$, consider $h:A\longrightarrow RLB$ and the associated 
universal function
$\zeta^{RLB}(A,LB):\textrm{Hom}_{\mathcal{C}}(A, RLB)\longrightarrow \textrm{Hom}_{\mathcal{D}}(LA, LB)$ 
in order to obtain

\begin{equation}\label{2311060020}
h^{RLB} := R\zeta^{RL}(A,LB)(h) = R\zeta^{RLB}(h)
\end{equation}

$\phantom{br}$

Therefore, the definition of the image of the 2-functor on 0-cells is

\begin{equation}\label{2311060027}
\Phi(\mathcal{C}, R, \eta^{RL}, \zeta^{R}(\sim)) := (\mathcal{C}, RL, \eta^{RL}, 
(\sim)^{RL}:= R\zeta^{RL}(\sim))
\end{equation}

$\phantom{br}$

\textbf{Note}: For a universal arrow $(\mathcal{D}, R^{\dagger},\eta^{RL\dagger})$ and a morphism 
$v:D\longrightarrow R^{\dagger}E$, the universal associated function 
$\zeta^{RL\dagger}(D,E)(v)$ is short denoted as $\zeta^{R\dagger E}(v)$ and for a morphism 
$w:D\longrightarrow R^{\dagger}L^{\dagger}\!E$ the function $\zeta^{RL\dagger}(D,L^{\dagger}\!E)(w)$
is short denoted as $\zeta^{RL\dagger E}(w)$. As before, in order to form the extensive induced
function $w^{RL\dagger E} := R^{\dagger}\zeta^{RL\dagger E}(w)$.\\

It has to be proved that the 2-functor is well-defined on 0-cells. That is to say,
$(\mathcal{C}, RL, \eta^{RL}, (\sim)^{RL})$ is a monad in extensive form, where
$(\sim)^{RL}:= R\zeta^{RL}(\sim)$.\\

\newtheorem{2311080057}{Proposition}[subsection]
\begin{2311080057}\label{2311080057}
$(\mathcal{C}, RL, \eta^{RL}, (\sim)^{RL})$ is an extensive monad.
\end{2311080057}

\noindent \emph{Proof}:\\

It only remains to prove the properties for the induced function 
$(\sim)^{RLB}:\textrm{Hom}_{\mathcal{C}}(A, RLB)\longrightarrow \textrm{Hom}_{\mathcal{C}}(RLA, RLB)$
and these are fulfilled as follows.\\

\begin{enumerate}

\item [i)] If the morphism for the right-hand side of the following universal triangle is 
$\eta^{RL}A:A\longrightarrow RLA$ then

\begin{equation*}
\xy<1cm,1cm>
\POS (16,24) *+{A} = "a12",
\POS (0, 0) *+{RLA} = "a21",
\POS (32, 0) *+{RLA} = "a23", 
\POS "a12" \ar_{\eta^{RL}\!A} "a21",
\POS "a12" \ar^{\eta^{RL}\!A} "a23",
\POS "a21" \ar_{R\zeta^{RL\!A}(\eta^{RL}\!A)} "a23"
\endxy
\end{equation*}

\noindent but there is also another morphism that makes the triangle commuting, 
namely $R1_{LA}$. Therefore, 

\begin{equation*}
(\eta^{RL}A)^{RLA} = R\zeta^{RLA}(\eta^{RL}A) = R(1_{LA}) = 1_{RLA}
\end{equation*}

\item [ii)] Consider $h:A\longrightarrow RLB$ as the right-hand side of the following 
universal triangle,

\begin{equation*}
\xy<1cm,1cm>
\POS (16,24) *+{A} = "a12",
\POS (0, 0) *+{RLA} = "a21",
\POS (32, 0) *+{RLB} = "a23", 
\POS "a12" \ar_{\eta^{RL}\!A} "a21",
\POS "a12" \ar^{h} "a23",
\POS "a21" \ar_{R\zeta^{RL\!B}(h)} "a23"
\endxy
\end{equation*}

Due to the commutativity of the triangle 
$h^{RLB}\!\cdot \eta^{RL}\!A = R\zeta^{RLB}(h)\cdot\eta^{RL}A = h$, which is the equation 
looked after.\\
 
\textbf{Note}: If a morphism of type $v:A\longrightarrow RX$ were considered instead, 
the commutativity of the triangle would provide the equation

\begin{equation}\label{2312110002}
R\zeta^{RX}(v)\cdot\eta^{RL}\!A = v
\end{equation}

$\phantom{br}$

\item [iii)] For $h:A\longrightarrow RLB$ and $k:B\longrightarrow RLC$, 
the equation $(k^{RLC}\!\cdot h)^{RLC} = k^{RLC}\!\cdot h^{RLB}$ has to be fulfilled. 
Without losing any generality, it will be considered the morphism $k:B\longrightarrow RX$ 
and the family of functions $\zeta^{RX}$ such that 
$(k^{RX}\!\cdot h)^{RX} = k^{RX}\!\cdot h^{RLB}$.\\

\begin{equation}\label{2312122247}
\begin{array}{ccc}
\xy
\POS (22,32) *+{A} = "a12",
\POS (33,16) *+{RLB} = "a23",
\POS (0, 0) *+{RLC} = "a31",
\POS (44, 0) *+{RX} = "a34",
\POS "a12" \ar_{\eta^{RL}C} "a31",
\POS "a31" \ar_{R\zeta^{RX}(k^{RX}\!\cdot h)} "a34",
\POS "a31" \ar^{h^{RL\!B}} "a23",
\POS "a12" \ar^{h} "a23",
\POS "a23" \ar^{k^{RX}} "a34",
\endxy
&, &
\xy
\POS (20,24) *+{A} = "a12",
\POS (0, 0) *+{RLC} = "a31",
\POS (40, 0) *+{RLB} = "a33",
\POS (70,0) *+{RLC} = "a34",
\POS "a12" \ar_{\eta^{RL}C} "a31",
\POS "a31" \ar^{h^{RL\!B}}_{R\zeta^{RL\!B}(h)} "a33",
\POS "a33" \ar^{k^{RX}}_{R\zeta^{RX}(k)} "a34",
\POS "a12" \ar^{h} "a33",
\endxy
\end{array}
\end{equation}

It is observed that the small universal triangle along with the $k^{RX}$ pull-forward 
composition can be inserted inside the big universal triangle making both commuting. 
Therefore both $R$-bases are equal which means that 
$R\zeta^{RX}(k^{RX}\!\cdot h) = R\zeta^{RX}(k)\cdot R\zeta^{RL\!B}(h)$ 
or using the corresponding definitions

\begin{equation}\label{2312110102}
(k^{RX}\!\cdot h)^{RX} = k^{RX}\!\cdot h^{RLB}
\end{equation}

Making $X = LC$, it is obtained instead

\begin{flushleft}
\color{lightgray}
\textbf{erasable} the following element has label 2312110103
\end{flushleft}

\begin{equation}\label{2312110103}
(k^{RLC}\!\cdot h)^{RLC} = k^{RLC}\!\cdot h^{RLB}
\end{equation}

\end{enumerate}

\begin{flushright}
$\Box$
\end{flushright}

\textbf{Note}: In the previous diagram, some arrows were labeled twice. This would mean that it 
is important to take into account, for the time being, both expressions of the same arrow. 

\subsection{On morphisms}\label{2311110913}

Let us define the 2-functor on 1-cells or morphisms, according to the previous
section, this definition has to have the following type.

\begin{equation}\label{2311111207}
\Phi(J, V,\rho^{JV}) := (J, (\sim)^{JRL}):
(\mathcal{C}, (\sim)^{RL})\longrightarrow (\mathcal{D}, (\sim)^{RL\dagger})
\end{equation}

$\phantom{br}$

\noindent and the corresponding extensive function is defined as follows.
For $p:D\longrightarrow JRLA$ construct the following diagram

\begin{equation*}
\xy<1cm,1cm>
\POS (22,32) *+{D} = "a12",
\POS (33,16) *+{JRLA} = "a23",
\POS (0, 0) *+{R^{\dagger}L^{\dagger}D} = "a31",
\POS (44, 0) *+{R^{\dagger}V\!LA} = "a34",
\POS (75, 0) *+{JRLA} = "a35",
\POS "a12" \ar_{\eta^{RL\dagger}D} "a31",
\POS "a12" \ar^{p} "a23",
\POS "a23" \ar^{\rho^{JV}\!LA} "a34",
\POS "a31" \ar_{R^{\dagger}\zeta^{R\dagger V\!L\!A}(\rho^{JV}\!LA\cdot p)} "a34",
\POS "a34" \ar_{\varrho^{JV}\!LA} "a35",
\POS "a31" \ar@/_2.5pc/_{p^{JRLA}:=} "a35",
\endxy
\end{equation*}

Therefore, 
$p^{JRLA}:= \varrho^{JV}\!LA\cdot R^{\dagger}\zeta^{R\dagger V\!L\!A}(\rho^{JV}\!LA\cdot p)$. 
Obviously, $\rho^{JV}\!LA\cdot p^{JRLA} = R^{\dagger}\zeta^{R\dagger V\!L\!A}(\rho^{JV}\!LA\cdot p)$.\\

\textbf{Note}: Only the triangle is universal, if it is pull-forward by the morphism 
$\varrho^{JV}\!LA$ it might not fulfill this property anymore.\\

The following pair of Propositions will show that the 2-functor is well-defined on 
1-cells.

\newtheorem{2311111214}{Proposition}[subsection]
\begin{2311111214}\label{2311111214}
$(JRLA, (\sim)^{JRL})$ is a $T$-algebra.
\end{2311111214}

\begin{flushleft}
\emph{Proof}:
\end{flushleft}

\begin{enumerate}

\item [i)] For $p:D\longrightarrow JRLA$, the condition $p^{JRLA}\cdot\eta^{RL\dagger}D = p$ 
is fulfilled by definition.

\item [ii)] The condition $p^{JRLA}\cdot r^{RL\dagger D} = (p^{JRLA}\cdot r)^{JRLA}$, for 
$p:D\longrightarrow JRLA$ and $r:Z\longrightarrow R^{\dagger}L^{\dagger}D$, is fulfilled 
as shown by the following diagram.

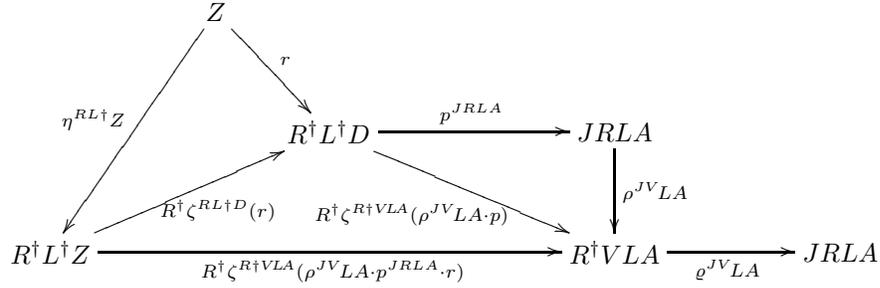
\begin{equation*}
\xy<1cm,1cm>
\POS (22,32) *+{Z} = "a12",
\POS (37,16) *+{R^{\dagger}L^{\dagger}D} = "a23",
\POS (0, 0) *+{R^{\dagger}L^{\dagger}Z} = "a31",
\POS (75, 16) *+{JRLA} = "a34",
\POS (75, 0) *+{R^{\dagger}VLA} = "a35",
\POS (105, 0) *+{JRLA} = "a36",
\POS (48,5) *+{\scriptstyle{R^{\dagger}\zeta^{R\dagger V\!L\!A}(\rho^{JV}\!LA\cdot p)}} = "x",
\POS "a12" \ar_{\eta^{RL\dagger}Z} "a31",
\POS "a31" \ar_{R^{\dagger}\zeta^{R\dagger V\!L\!A}(\rho^{JV}\!LA\cdot p^{JRLA}\cdot r)} "a35",
\POS "a31" \ar_{R^{\dagger}\zeta^{RL\dagger D}(r)} "a23",
\POS "a23" \ar "a35",
\POS "a12" \ar^{r} "a23",
\POS "a23" \ar^{p^{JRLA}} "a34",
\POS "a34" \ar^{\rho^{JV}\!LA} "a35",
\POS "a35" \ar_{\varrho^{JV}\!LA} "a36",
\endxy
\end{equation*}

$\phantom{br}$

In the same way as before \eqref{2312122247}, the small universal triangle with right-hand side $r$ 
can be inserted into the universal triangle with right-hand side $\rho^{JV}\!LA\cdot p^{JRLA}\!\cdot r$
and the resulting parallel $R^{\dagger}$-bases make the big universal triangle commuting. Therefore
pulling forward with the morphism $\varrho^{JV}\!LA$ gives

\begin{eqnarray*}
(p^{JRLA}\cdot r)^{JRLA} 
&=& \varrho^{JV}\!LA\cdot R^{\dagger}\zeta^{R\dagger V\!L\!A}(\rho^{JV}\!LA\cdot p^{JRLA}\cdot r)\\
&=& \varrho^{JV}\!LA\cdot R^{\dagger}\zeta^{R\dagger V\!L\!A}(\rho^{JV}\!LA\cdot p)\cdot R^{\dagger}\zeta^{RL\dagger D}(r)\\
&=& p^{JRLA}\cdot R^{\dagger}\zeta^{RL\dagger D}(r)\\
&=& p^{JRLA}\cdot r^{RL\dagger D}
\end{eqnarray*}

\end{enumerate}

\begin{flushright}
$\Box$
\end{flushright} 

\newtheorem{2311111219}[2311111214]{Proposition}
\begin{2311111219}\label{2311111219}
$Jh^{JRLB}$ is a $T$-algebra morphism for any $h:A\longrightarrow RLB$.
\end{2311111219}

\begin{flushleft}
\emph{Proof}:
\end{flushleft}

Consider a morphism of type $p:D\longrightarrow JRLA$ and the following diagram 

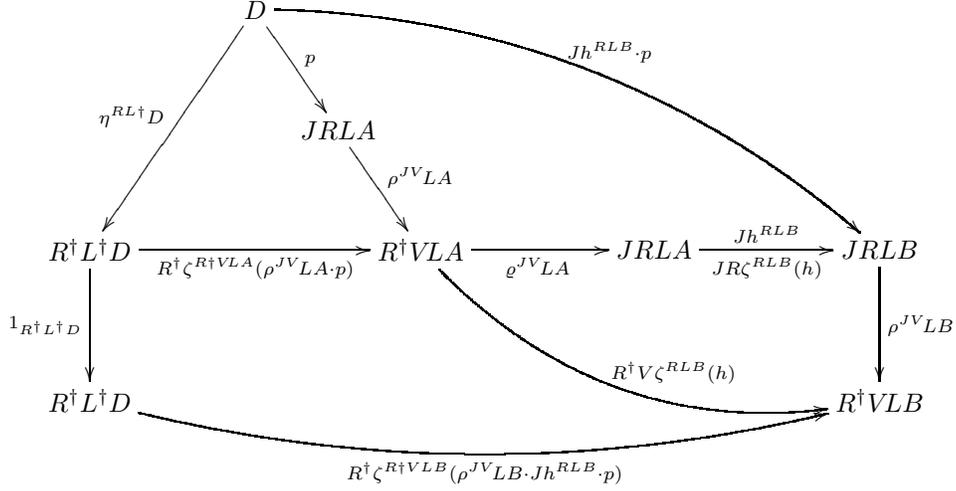
\begin{equation*}
\xy<1cm,1cm>
\POS (22,32) *+{D} = "a12",
\POS (33,16) *+{JRLA} = "a23",
\POS (0, 0) *+{R^{\dagger}L^{\dagger}D} = "a31",
\POS (44, 0) *+{R^{\dagger}V\!LA} = "a34",
\POS (75, 0) *+{JRLA} = "a35",
\POS (105,0) *+{JRLB} = "a36",
\POS (0,-20) *+{R^{\dagger}L^{\dagger}D} = "a41",
\POS (105,-20) *+{R^{\dagger}V\!LB} = "a46",
\POS "a12" \ar_{\eta^{RL\dagger}D} "a31",
\POS "a31" \ar_{1_{R^{\dagger}L^{\dagger}D}} "a41",
\POS "a41" \ar@/_1.7pc/_{R^{\dagger}\zeta^{R\dagger VLB}(\rho^{JV}\!LB\cdot Jh^{RLB}\cdot p)} "a46",
\POS "a31" \ar_{R^{\dagger}\zeta^{R\dagger V\!LA}(\rho^{JV}\!LA\cdot p)} "a34",
\POS "a34" \ar@/_2.2pc/^{\phantom{012}R^{\dagger}V\zeta^{RLB}(h)} "a46",
\POS "a34" \ar_{\varrho^{JV}\!LA} "a35",
\POS "a35" \ar_{JR\zeta^{RLB}(h)}^{Jh^{RLB}} "a36",
\POS "a36" \ar^{\rho^{JV}\!LB} "a46",
\POS "a12" \ar^{p} "a23",
\POS "a23" \ar^{\rho^{JV}\!LA} "a34",
\POS "a12" \ar@/^2pc/^{Jh^{RLB}\cdot p} "a36",
\endxy
\end{equation*}

$\phantom{br}$

As before, the small universal triangle with right-hand side $\rho^{JV}\!LA\cdot p$ and pulled 
forward by the morphism $\rho^{JV}\!LB\cdot JR\zeta^{RLB}(h)\cdot \varrho^{JV}\!LA$ can be inserted 
into the universal triangle with right-hand side $\rho^{JV}\!LB\cdot Jh^{RLB}\!\cdot p$.\\

On the other hand, by naturality of $\rho^{JV}\!L$ over the morphism $\zeta^{RLB}(h)$, 
the following equation holds 
$R^{\dagger}V\zeta^{RLB}(h)\cdot \rho^{JV}\!LA = \rho^{JV}\!LB\cdot JR\zeta^{RLB}(h)$. 
Using this equality, the exterior  universal triangle with right-hand side 
$\rho^{JV}\!LB\cdot Jh^{RLB}\!\cdot p$ 
has a pair of parallel $R^{\dagger}$-bases which make this universal triangle commute. Therefore 
this pair gives the equality  
$R^{\dagger}\zeta^{R\dagger VLB}(\rho^{JV}\!LB\cdot Jh^{RLB}\cdot p) = 
R^{\dagger}V\zeta^{RLB}(h)\cdot R^{\dagger}\zeta^{R\dagger VLA}(\rho^{JV}\!LA\cdot p)$ 
and if this equality is composed, in turn, with $\varrho^{JV}\!LB$ then \\

\begin{eqnarray*}
(Jh^{RLB}\cdot p)^{JRLB} 
&=& \varrho^{JV}\!LB\cdot R^{\dagger}\zeta^{R\dagger VLB}(\rho^{JV}\!LB\cdot Jh^{RLB}\cdot p)\\
&=& \varrho^{JV}\!LB\cdot R^{\dagger}V\zeta^{RLB}(h)\cdot R^{\dagger}\zeta^{\dagger VLA}(\rho^{JV}\!LA\cdot p)\\
&=& JR\zeta^{RLB}(h)\cdot \varrho^{JV}\!LA\cdot R^{\dagger}\zeta^{R\dagger VLA}(\rho^{JV}\!LA\cdot p)\\
&=& Jh^{RLB}\cdot p^{JRLA}
\end{eqnarray*}

$\phantom{br}$

\noindent where in the third equality, the naturality of $\rho^{JV}$ over the morphism 
$\zeta^{RLB}(h)$ was used. 

\begin{flushright}
$\Box$
\end{flushright}

\subsection{On transformations of universal arrows}\label{2311142003}

The image of the 2-functor $\Phi$, on 2-cells, is defined as follows

\begin{equation}\label{2311142007}
\Phi(\alpha, \beta) := \alpha:(J, (\sim)^{JRL})\longrightarrow (K, (\sim)^{KRL})
\end{equation}

$\phantom{br}$

Let us prove that it fulfills the condition for a 2-cell in $\mathbf{EMnd}(Cat)$ 
through the following proposition.

\newtheorem{2311142012}{Proposition}[subsection]
\begin{2311142012}\label{2311142012}
The component $\alpha RLA: (JRLA, (\sim)^{JRLA})\longrightarrow (KRLA, (\sim)^{KRLA}$ is a 
morphism of $RL^{\dagger}$-algebras.
\end{2311142012}

\begin{flushleft}
\emph{Proof}:
\end{flushleft}

Consider a morphism of type $p:D\longrightarrow JRLA$, it has to be proved that 
$\alpha RLA\cdot p^{JRLA} = (\alpha RLA\cdot p)^{KRLA}$. In order to do so, consider 
the following universal triangle with right-hand side $\rho^{KW}\!LA\cdot\alpha RLA\cdot p$

\begin{equation*}
\xy<1cm,1cm>
\POS (22,32) *+{D} = "a12",
\POS (33,16) *+{JRLA} = "a23",
\POS (62,16) *+{KRLA} = "a25",
\POS (0, 0) *+{R^{\dagger}L^{\dagger}D} = "a31",
\POS (44, 0) *+{R^{\dagger}V\!LA} = "a34",
\POS (75, 0) *+{R^{\dagger}WLA} = "a36",
\POS "a12" \ar_{\eta^{RL\dagger}D} "a31",
\POS "a31" \ar@/_2.5pc/_{R^{\dagger}\zeta^{R\dagger W\!L\!A}(\rho^{KW}\!LA\cdot\alpha RLA\cdot p)} "a36",
\POS "a12" \ar^{p} "a23",
\POS "a23" \ar^{\rho^{JV}\!LA} "a34",
\POS "a31" \ar_{R^{\dagger}\zeta^{R\dagger V\!L\!A}(\rho^{JV}\!LA\cdot p)} "a34",
\POS "a34" \ar_{R^{\dagger}\beta LA} "a36",
\POS "a23" \ar^{\alpha RLA} "a25",
\POS "a25" \ar^{\rho^{KW}\!LA} "a36",
\endxy
\end{equation*}

$\phantom{br}$

Then, 
$R^{\dagger}\beta LA\cdot R^{\dagger}\zeta^{R\dagger V\!L\!A}(\rho^{JV}\!LA\cdot p) 
= R^{\dagger}\zeta^{R\dagger W\!LA}(\rho^{KW}\!LA\cdot\alpha RLA\cdot p)$ since they both
make the corresponding universal triangle commute at the same time with a 
parallel $R^{\dagger}$-base. Note that the compatilibity condition required for the
pair $(\alpha,\beta)$, as in \eqref{2311080011}, was used. Furthermore, if this equality 
is pull forward by $\varrho^{KW}\!LA$ then 

\begin{eqnarray*}
(\alpha RLA\cdot p)^{KRLA} 
&=& \varrho^{KW}\!LA\cdot R^{\dagger}\zeta^{\dagger WLA}(\rho^{KW}\!LA\cdot\alpha RLA\cdot p)\\
&=& \varrho^{KW}\!LA\cdot R^{\dagger}\beta LA\cdot R^{\dagger}\zeta^{R\dagger V\!LA}(\rho^{JV}\!LA\cdot p)\\
&=& \alpha RLA\cdot \varrho^{JV}\!LA\cdot R^{\dagger}\zeta^{R\dagger V\!LA}(\rho^{JV}\!LA\cdot p)\\
&=& \alpha RLA\cdot p^{JRLA}
\end{eqnarray*}

In the third equality the compatibility condition for pair $(\alpha, \beta)$ was also applied. 
The rest of the equalities only make use of the definition of their respective functions, 
for example $(\sim)^{KRL}$.

\begin{flushright}
$\Box$
\end{flushright}

\newtheorem{2312292248}[2311142012]{Proposition}
\begin{2312292248}\label{2312292248}
There exists a 2-functor of type $\Phi:\mathbf{UArr}(Cat)\longrightarrow\mathbf{EMnd}(Cat)$.
\end{2312292248}

\begin{flushleft}
\emph{Proof}:
\end{flushleft}

Let us prove the functoriality on composition 1-cells. Take for instance 
$(J,V,\rho^{JV}):(\mathcal{C},R)\longrightarrow(\mathcal{D},R^{\dagger})$ and
$(K,W,\rho^{KW}):(\mathcal{D},R^{\dagger})\longrightarrow(\mathcal{E},R^{\ddagger})$, 
it has to be proved that 
$\Phi(K,W,\rho^{KW})\cdot\Phi(J,V,\rho^{JV}) = \Phi(JK,VW,\rho^{KW}V\circ K\rho^{JV})$. 
The following diagram shows this equality for $x:E\longrightarrow KJRLA$

\begin{equation}
\xy
\POS (0,50) *+{E} = "a11",
\POS (60,50) *+{KJRLA} = "a12",
\POS (60,30) *+{KR^{\dagger}L^{\dagger}JRLA} = "a22",
\POS (120,30) *+{KR^{\dagger}VLA} = "a23",
\POS (0,10) *+{R^{\ddagger}L^{\ddagger}E} = "a31",
\POS (60,10) *+{R^{\ddagger}WL^{\dagger}JRLA} = "a32",
\POS (120,10) *+{R^{\ddagger}WVLA} = "a33",
\POS (0,0) *+{R^{\ddagger}L^{\ddagger}E} = "a41",
\POS (120,0) *+{R^{\ddagger}WVLA} = "a43",
\POS (30,6) *+{\scriptstyle{R^{\ddagger}\zeta^{R\ddagger W\!JRLA}(\rho^{KW}\!L^{\dagger}JRLA\cdot K\eta^{RL\dagger}\!JRLA\cdot x)}} = "x",
\POS "a11" \ar_{\eta^{RL\ddagger}E} "a31",
\POS "a31" \ar_{1_{R^{\ddagger}L^{\ddagger}E}} "a41",
\POS "a41" \ar_{R^{\ddagger}\zeta^{R\ddagger WLA}(\rho^{KW}\!V\!LA\cdot K\rho^{JV}\!LA\cdot x)} "a43",
\POS "a31" \ar "a32",
\POS "a32" \ar_{R^{\ddagger}W\zeta^{R\dagger VLA}(\rho^{JV}\!LA\cdot 1_{JRLA})} "a33",
\POS "a11" \ar^{x} "a12",
\POS "a12" \ar^{K\eta^{RL\dagger}JRLA} "a22",
\POS "a22" \ar@<2.5pt>^{\rho^{KW}\!L^{\dagger}\!JRLA} "a32",
\POS "a32" \ar@<2.5pt>^{\varrho^{KW}\!L^{\dagger}\!JRLA} "a22",
\POS "a22" \ar_{KR^{\dagger}\zeta^{R\dagger V\!LA}(\rho^{JV}\!LA\cdot 1_{JRLA})} "a23",
\POS "a23" \ar^{\rho^{KW}\!V\!LA} "a33",
\POS "a12" \ar@/^1.2pc/^{\phantom{123}K\rho^{JV}\!LA} "a23",
\POS "a33" \ar^{1_{R^{\ddagger}WV\!LA}} "a43"
\endxy
\end{equation}

\noindent where the compositions of images corresponds, first to the $R^{\ddagger}$-base of the left
 rectangle (universal triangle) and then to the $KR^{\dagger}$-base of the small $K$ universal upper triangle. 
By appling naturality of $\rho^{KW}$ over the morphism $\zeta^{R\dagger VLA}(\rho^{JV}\!LA\cdot 1_{JRLA})$
it is obtained a pair of parallel $R^{\ddagger}$-bases for the big universal triangle with right-hand
side $\rho^{KW}VLA\cdot K\rho^{JV}\!LA\cdot x$. This will lead to the required equality.\\

Since the left 2-functor is the identity on 2-cells, there remains nothing to check.

\begin{flushright}
$\Box$
\end{flushright}

\section{The Right 2-functor $\Psi:\mathbf{EMnd}(Cat)\longrightarrow\mathbf{UArr}(Cat)$}\label{2312151908}

In this section, the right 2-functor $\Psi:\mathbf{EMnd}(Cat)\longrightarrow\mathbf{UArr}(Cat)$ 
is constructed on each $n$-cell. Remember that each subsection corresponds to a specific $n$-cell.

\subsection{On extensive monads}\label{2311060048}

Let us construct the image of an extensive monad for the right 2-functor $\Psi$, that is to say a 
universal arrow.\\

The proposal for the universal arrow, at $A$ in $\mathcal{C}$, has the following type 
$\eta^{UFS}\!A:A\longrightarrow U^{S}F^{S}\!A$. Since $U^{S}F^{S}\!A = U^{S}(SA, (\sim)^{S\! A}) = SA$ then 
the original morphism $\eta^{S}\!A:A\longrightarrow SA$ can be proposed as the universal arrow.\\

In order to have $\eta^{S}A$ as a universal arrow from $A$ to $U^{S}$ then, for any 
$v:A\longrightarrow U^{S}(N, (\sim)^{N})$, there must exist a unique morphism of $S$-algebras
$w:F^{S}A = (SA,(\sim)^{SA})\longrightarrow (N, (\sim)^{N})$ such that the following triangle commutes

\begin{equation}\label{2312180021}
\xy
\POS (20,27) *+{A} = "a12",
\POS (0, 0) *+{U^{S}F^{S}\!A} = "a21",
\POS (40, 0) *+{U^{S}(N, (\sim)^{N})} = "a23",
\POS "a12" \ar_{\eta^{S}\!A} "a21",
\POS "a12" \ar^{v} "a23",
\POS "a21" \ar_-{U^{S}w} "a23"
\endxy
\end{equation}\\

Consider the morphism $v:A\longrightarrow N = U^{S}(N, (\sim)^{N})$ and from the properties for the algebra
$(N, (\sim)^{N})$ there exist, through the associated function of the algebra, 
a morphism of type $v^{N}:SA\longrightarrow N$. It is proved that it is a morphism of $S$-algebras.\\

\newtheorem{2311111232}{Proposition}[subsection]
\begin{2311111232}
The induced morphism $v^{N}:SA\longrightarrow N$, is a morphism of algebras 
$v^{N}:(SA, (\sim)^{SA})\longrightarrow (N, (\sim)^{N})$.
\end{2311111232}

\begin{flushleft}
\emph{Proof}:
\end{flushleft}

Let $h:X\longrightarrow SA$ be an arbitrary morphism, then it has to be proved that 
$v^{N}\cdot h^{SA} = (v^{N}\cdot h)^{N}$ but this is one of the conditions for $(N, (\sim)^{N})$ 
to be an algebra and, in particular, for $v$ to be an algebraic morphism, check Subsection 
\eqref{2311141948}.

\begin{flushright}
$\Box$
\end{flushright}

The required morphism $w$ is then $v^{N}$. The condition $v^{N}\!\cdot\eta^{S}A = v$ is fulfilled 
from the properties of the algebra $(N,(\sim)^{N})$ and the universal triangle commutes. 
The uniqueness of $w$ is proved as follows. \\

Consider another morphism of $S$-algebras
$\widetilde{w}:(SA, (\sim)^{S})\longrightarrow (N, (\sim)^{N})$ such that 
$\widetilde{w}\cdot\eta^{S}A = v$ , then 

\begin{equation*}
\widetilde{w} 
= \widetilde{w}\cdot 1_{SA} 
= \widetilde{w}\cdot(\eta^{S}A)^{SA} 
= (\widetilde{w}\cdot\eta^{S}A)^{N} = v^{N},\\
\end{equation*}

\noindent where the third equality follows since $\widetilde{w}$ is a morphism of 
$S$-algebras.\\

Finally, the universal function 
$\zeta^{U\!S}(A,(N, (\sim)^{N})):\mathrm{Hom}_{\mathcal{C}}(A,N)\longrightarrow
\mathrm{Hom}_{\mathcal{C}^{S}}((SA, (\sim)^{SA}),(N, (\sim)^{N}))$ 
is defined, using the respective short notation, $\zeta^{USN}(v) = v^{N}$. 
This completes the structure for a universal arrow.

\begin{equation}\label{2312152122}
\Psi(\mathcal{C},S,\eta^{S},(\sim)^{S}) := (\mathcal{C}, U^{S}, \eta^{U\!FS}, \zeta^{U\!S}(\sim))
\end{equation}

\noindent where $\eta^{UFS}:= \eta^{S}$.

\subsection{On morphisms of extensive monads}\label{2311162317}

Consider a morphism of extensive monads 
$(P, (\sim)^{PS}): (\mathcal{C}, (\sim)^{S})\longrightarrow(\mathcal{D}, (\sim)^{T})$ 
then a morphism of univeral arrows of the following type
$(P, \widehat{P}, \rho^{P}): (\mathcal{C},U^{S}, \eta^{U\!FS}, \zeta^{U\!S})
\longrightarrow (\mathcal{D},U^{T}, \eta^{U\!FT}, \zeta^{U\!T})$
has to be constructed.\\

First of all, the following functor is defined 
$\widehat{P}:\mathcal{C}^{S}\longrightarrow\mathcal{D}^{T}$. 

\begin{enumerate}

\item [\emph{i})] For objects, consider $(N,(\sim)^{N})$ 

\begin{equation}\label{2311162331}
\widehat{P}(N, (\sim)^{N}) := (PN, (\sim)^{PN})
\end{equation}

In order to define the respective function, consider the morphism $q:D\longrightarrow PN$ 
and use the structure of the algebra $(PSN, (\sim)^{PSN})$ to get 
$(P\eta^{S}N\!\cdot q)^{PSN}\!:TD\longrightarrow PSN$ and finally do the pull forward 
with $P1_{N}^{\:N}$ 

\begin{equation}\label{2311170121}
q^{PN}\!:= P1_{N}^{\:N}\cdot(P\eta^{S}N\cdot q)^{PSN}:TD\longrightarrow PN
\end{equation}

$\phantom{br}$

\item [\emph{ii})] On morphisms, consider $m:(M,(\sim)^{M})\longrightarrow(N,(\sim)^{N})$ 
then 

\begin{equation}\label{2311170124}
\widehat{P}m := Pm: (PM, (\sim)^{PM})\longrightarrow (PN, (\sim)^{PN})
\end{equation}

\end{enumerate}

\newtheorem{2311170127}{Proposition}[subsection]
\begin{2311170127}\label{2311170127}
$\widehat{P}:\mathcal{C}^{S}\longrightarrow\mathcal{D}^{T}$ is a functor.
\end{2311170127}

\begin{flushleft}
\emph{Proof}:
\end{flushleft}

First $(PN, (\sim)^{PN})$ is a $T$-algebra, then 

\begin{enumerate}

\item [\emph{i})] Consider $q:D\longrightarrow PN$, 
$P1_{N}^{\:N}\cdot(P\eta^{S}N\cdot q)^{PSN}\!\!\cdot\eta^{T}D = P1_{N}^{\:N}\cdot P\eta^{S}N\cdot q = q$.
The equality holds because $PSN$ is a $T$-algebra and $N$ is a $S$-algebra. 

\item [\emph{ii})] For $r:Z\longrightarrow TD$ and $q:D\longrightarrow PN$.

\begin{eqnarray*}
(q^{PN}\cdot r)^{PN} 
&=& P1_{N}^{\:N}\cdot[P\eta^{S}N\cdot P1_{N}^{\:N}\cdot(P\eta^{S}N\cdot q)^{PSN}\cdot r]^{PSN}\\  
&=& P1_{N}^{\:N}\cdot[P(\eta^{S}N\cdot 1_{N}^{\:N})^{SN}\cdot P\eta^{S}N\cdot(P\eta^{S}N\cdot q)^{PSN}\cdot r]^{PSN}\\
&=& P1_{N}^{\:N}\cdot P(\eta^{S}N\cdot 1_{N}^{\:N})^{SN}\cdot[P\eta^{S}N\cdot(P\eta^{S}N\cdot q)^{PSN}\cdot r]^{PSN}\\
&=& P(1_{N}^{\:N}\cdot\eta^{S}N\cdot 1_{N}^{\:N})^{N}\cdot[P\eta^{S}N\cdot(P\eta^{S}N\cdot q)^{PSN}\cdot r]^{PSN}\\
&=& P(1_{N}^{\:N})^{N}\cdot[P\eta^{S}N\cdot(P\eta^{S}N\cdot q)^{PSN}\cdot r]^{PSN}\\
&=& P1_{N}^{\:N}\cdot P1_{SN}^{\:SN}\cdot[P\eta^{S}N\cdot(P\eta^{S}N\cdot q)^{PSN}\cdot r]^{PSN}\\
&=& P1_{N}^{\:N}\cdot [P1_{SN}^{\:SN}\cdot P\eta^{S}N\cdot(P\eta^{S}N\cdot q)^{PSN}\cdot r]^{PSN}\\
&=& P1_{N}^{\:N}\cdot [(P\eta^{S}N\cdot q)^{PSN}\cdot r]^{PSN}\\
&=& P1_{N}^{N\:}\cdot (P\eta^{S}N\cdot q)^{PSN}\cdot r^{TD}\\
&=& q^{PN}\cdot r^{TD}
\end{eqnarray*}

The second equality takes place since the property of type $h^{SB}\cdot \eta^{S}A = h$ was applied.
The third equality follows since $Ph^{SB}$ is a morphism of $T$-algebras for any $h:A\longrightarrow SB$. 
In the fourth equality, that $1_{N}$ is an algebraic morphism was used. 
In the fifth equality, the property $1_{N}^{N}\cdot\eta^{S}\!N = 1_{N}$ was applied.
For the sixth equality, that $1_{N}$ is an algebraic morphism was used again over the morphism $1_{SN}$.
The seventh equality follows since $P1_{SN}^{\:SN}$ is morphism of $T$-algebras. In the eight equality,
the property $h^{SB}\cdot \eta^{S}A = h$ was applied again. Finally, in the nineth equality, that
$(PSN, (\sim)^{PSN})$ is a $T$-algebra was applied.\\

Second, in order to prove that $Pm: (PM, (\sim)^{PM})\longrightarrow (PN, (\sim)^{PN})$ is a morphism
of $T$-algebras, consider $Pm: PM\longrightarrow PN$ and $q:D\longrightarrow PM$.

\begin{eqnarray*}
Pm\cdot q^{PM} 
&=& Pm\cdot P1_{M}^{\:M}\cdot(P\eta^{S}M\cdot q)^{PSM}\\
&=& P(m\cdot 1_{M})^{N}\cdot(P\eta^{S}M\cdot q)^{PSM}\\ 
&=& P(1_{N}\cdot m)^{N}\cdot(P\eta^{S}M\cdot q)^{PSM}\\
&=& P(1_{N}^{\:N}\cdot\eta^{S}N\cdot m)^{N}\cdot(P\eta^{S}M\cdot q)^{PSM}\\ 
&=& P1_{N}^{\:N}\cdot P(\eta^{S}N\cdot m)^{SN}\cdot(P\eta^{S}M\cdot q)^{PSM}\\
&=& P1_{N}^{\:N}\cdot[P(\eta^{S}N\cdot m)^{SN}\cdot P\eta^{S}M\cdot q]^{PSN}\\ 
&=& P1_{N}^{\:N}\cdot[P\eta^{S}N\cdot Pm\cdot q]^{PSN}\\
&=& (Pm\cdot q)^{PN} 
\end{eqnarray*}

In the second equality, the property for $m$ to be a morphism of $S$-algebras was applied.
In the fourth equality, the property $1_{N}^{N}\cdot\eta^{S}N = 1_{N}$ was used. In the fifth equality,
it was used the property for $1_{N}$ being an algebraic morphism, over $N$. In the sixth equality,
the requirement for $P(\eta^{S}N\cdot m)^{SN}$ to be a morphism of $T$-algebras was applied. Finally,
in the seventh equality, the property $(\eta^{S}N\cdot m)^{SN}\cdot\eta^{S}M = \eta^{S}N\cdot m$ was
used.

\end{enumerate}

\begin{flushright}
$\Box$
\end{flushright}

The associated natural transformation has the following type 
$\rho^{P}:PU^{S}\longrightarrow U^{T}\widehat{P}:\mathcal{C}^{S}\longrightarrow\mathcal{D}$ 
, whose component at $(N,(\sim)^{N})$ is $\rho^{P}(N,(\sim)^{N}):PN\longrightarrow PN$. 
Because of the previous, the component is defined as the identity at $N$, 
$\rho^{P}(N,(\sim)^{N}) = 1_{N}$. This finalizes the definition for the image of the 
2-functor on 1-cells, then

\begin{equation}\label{2312161908}
\Psi(P, (\sim)^{PS}) 
:= (P, \widehat{P}, \rho^{P}): 
(\mathcal{C}, U^{S}, \eta^{UFS}, \zeta^{US}(\sim))\longrightarrow 
(\mathcal{D}, U^{T}, \eta^{UFT}, \zeta^{UT}(\sim))
\end{equation}

\subsection{On transformations of extensive monads}\label{2111202030}

Consider a transformation for extensive monads of the following type 
$\theta:(P,(\sim)^{P})\longrightarrow(Q,(\sim)^{Q}):
(\mathcal{C},(\sim)^{S})\longrightarrow(\mathcal{D},(\sim)^{T})$. 
Therefore define the following transformation 
$\widehat{\theta}:\widehat{P}\longrightarrow\widehat{Q}:
\mathcal{C}^{S}\longrightarrow\mathcal{D}^{T}$ 
through its components.

\begin{equation}\label{2311202111}
\widehat{\theta}(N, (\sim)^{N}) := \theta N: (PN, (\sim)^{PN})\longrightarrow(QN, (\sim)^{QN})
\end{equation}

$\phantom{br}$

It has to be checked that it is a morphism of $T$-algebras. In order to do this, consider 
$q:D\longrightarrow PN$ then it has to be proved that $\theta N\cdot q^{PN} = (\theta N\cdot q)^{QN}$.

\begin{eqnarray*}\label{2311202202}
\theta N\cdot q^{PN} 
&=& \theta N\cdot P1_{N}^{\:N}\cdot(P\eta^{S}N\cdot q)^{PSN}\\
&=& Q1_{N}^{\:N}\cdot\theta SN\cdot(P\eta^{S}N\cdot q)^{PSN}\\
&=& Q1_{N}^{\:N}\cdot(\theta SN\cdot P\eta^{S}N\cdot q)^{QSN}\\
&=& Q1_{N}^{\:N}\cdot(Q\eta^{S}N\cdot\theta N\cdot q)^{QSN}\\
&=& (\theta N\cdot q)^{QN}
\end{eqnarray*}

In the second equality, the naturality of $\theta$ over the morphism $P1_{N}^{N}$ was applied. 
In the third equality, it was used the fact that the component $\theta SN$ is a morphism of $T$-algebras. 
Finally, the naturality of $\theta$ over the morphism $\eta^{S}N$ was applied.\\

Therefore, 

\begin{equation}\label{2312161933}
\Psi(\theta) := (\theta, \widehat{\theta}): 
(P, \widehat{P}, \rho^{P})\longrightarrow(Q, \widehat{Q}, \rho^{Q}):
(\mathcal{C}, U^{S}, \eta^{UFS}, \zeta^{US}(\sim))\longrightarrow 
(\mathcal{D}, U^{T}, \eta^{UFT}, \zeta^{UT}(\sim))
\end{equation}

Now that the definition of the 2-functor on each $n$-cell has been given, it can be proved 
that fulfills the rest of the conditions. 

\newtheorem{2312161939}{Proposition}[subsection]
\begin{2312161939}\label{2312161939}
There is a 2-functor of type $\Psi:\mathbf{EMnd}(Cat)\longrightarrow\mathbf{UArr}(Cat)$.
\end{2312161939}

\begin{flushleft}
\emph{Proof}:
\end{flushleft}

Consider the vertical structure of the 2-category $\mathbf{EMnd}(Cat)$. That is to say, let there be
a pair of natural transformations of type $\theta:(P,(\sim)^{PS})\longrightarrow(X,(\sim)^{XS})$ and
$\vartheta:(X,(\sim)^{XS})\longrightarrow(Q,(\sim)^{QS})$. It was already checked that this vertical 
composition is well-defined. 

\begin{equation*}
\Psi(\vartheta\circ\theta)(N,(\sim)^{N})
= \widehat{(\vartheta\circ\theta)}(N,(\sim)^{N}) 
= (\vartheta\circ\theta)N 
= \vartheta N\cdot\theta N 
= \widehat{\vartheta}(N,(\sim)^{N})\cdot \widehat{\theta}(N,(\sim)^{N})
= [\Psi(\vartheta)\circ\Psi(\theta)](N,(\sim)^{N})
\end{equation*}

$\phantom{br}$

Before going into the whiskerings, let us check on the compositions of 1-cells. In order to do this,
consider $(P,(\sim)^{PS}):(\mathcal{C},S)\longrightarrow(\mathcal{D},T)$ and
$(W,(\sim)^{WT}):(\mathcal{D},T)\longrightarrow(\mathcal{X},U)$ then it is only important to consider
$\Psi(P,(\sim)^{PS}) = \widehat{P}$ and $\Psi(W,(\sim)^{WT}) = \widehat{W}$ and their composition 
$\Psi(WP,(\sim)^{WPS})$.\\

Remember that $\widehat{P}(N,(\sim)^{N}) = (PN,(\sim)^{PN})$, on the $S$-algebra $(N,(\sim)^{N})$ and for
a morphism of type $q:D\longrightarrow PN$, $q^{PN} = P1_{N}^{\:N}\cdot[P\eta^{S}N\cdot q]^{PSN}$. 
Likewise,  $\widehat{W}(I,(\sim)^{I}) = (WI,(\sim)^{WI})$, on the $T$-algebra $(I,(\sim)^{I})$ and
for a morphism of type $i:X\longrightarrow WI$, $i^{WI} = W1_{I}^{\:I}\cdot[W\eta^{T}I\cdot i]^{WTI}$.\\

Therefore if $I = PN$ and $i:X\longrightarrow WPN$

\begin{equation*}
i^{WPN} = W1_{PN}^{\:PN}\cdot[W\eta^{T}PN\cdot i]^{WTPN}
\end{equation*}

$\phantom{br}$

This corresponds to the compostion $(WP, (\sim)^{WPS}):(\mathcal{C},S)\longrightarrow(\mathcal{X},U)$
where the image of the associated function on $w:X\longrightarrow WPSA$ is  
$w^{WPSA} = W1_{PSA}^{\:\:PSA}\cdot[W\eta^{T}PSA\cdot w]^{WTPSA}$ \eqref{2312241907}.\\

Now, let us consider whiskerings. Take into account the following configuration in 
$\mathbf{EMnd}(Cat)$, that is to say a pull forward composition.\\ 

\begin{equation*}
\xy
\POS (0,0) *+{(\mathcal{C},S)} = "a11",
\POS (28,0) *+{(\mathcal{D},T)} = "a12",
\POS (58,0) *+{(\mathcal{X},U)} = "a13",
\POS (12,0) *+{\downarrow\:\theta} = "y",
\POS "a11" \ar@/^1.2pc/^{(P,(\sim)^{PS})} "a12",
\POS "a11" \ar@/_1.2pc/_{(Q,(\sim)^{QS})} "a12",
\POS "a12" \ar_{(W,(\sim)^{WT})} "a13",
\endxy
\end{equation*}

It is enough to check the following, since the corresponding 1-cells compositions
are equal.

\begin{equation*}
(\widehat{W}\widehat{\theta})(N,(\sim)^{N}) 
= \widehat{W}(\theta N) 
= W\theta N 
= \widehat{W\theta}(N,(\sim)^{N})
\end{equation*} 

$\phantom{br}$

In the previous equalities, the definitions \eqref{2311202111} and \eqref{2311170124} were applied.
The proof for the pull back composition is similar.

\begin{flushright}
$\Box$
\end{flushright}

\section{The Unit and the Counit of the 2-adjunction}\label{2311202228}

\subsection{The Unit}\label{2312172244}

Taking into account the previous pair of 2-functors, the unit of the 2-adjunction has 
the following type

\begin{equation}\label{2311202230}
\eta^{\Psi\Phi}: 
1_{\mathbf{UArr}}\longrightarrow \Psi\Phi:
\mathbf{UArr}(Cat)\longrightarrow\mathbf{UArr}(Cat)
\end{equation}

$\phantom{br}$

The components of the unit are defined as follows. Let $(\mathcal{C}, R, \eta^{RL},\zeta^{R}(\sim))$ be a 
universal arrow, where $R:\mathcal{X}\longrightarrow\mathcal{C}$. The image on this 
universal arrow, under the 2-functor $\Phi$, is the extensive monad
$\Phi(\mathcal{C}, R, \eta^{RL},\zeta^{RL}(\sim))=(\mathcal{C},RL,\eta^{RL},(\sim)^{RL})$. In turn,
the image of this extensive monad, under $\Psi$, is 
$\Psi(\mathcal{C},RL,\eta^{RL},(\sim)^{RL}) = (\mathcal{C},U^{RL},\eta^{URL},\zeta^{URL}(\sim))$. 
Therefore, the proposed component of the unit, at $(\mathcal{C}, R, \eta^{RL},\zeta^{R}(\sim))$, is

\begin{equation}
\eta^{\Psi\Phi E}(\mathcal{C}, R, \eta^{RL},\zeta^{R}(\sim)):=(1_{\mathcal{C}},K(\mathcal{C},R),\rho^{KR}):
(\mathcal{C}, R, \eta^{RL},\zeta^{R}(\sim))\longrightarrow(\mathcal{C},U^{RL},\eta^{URL},\zeta^{URL}(\sim)) 
\end{equation}

$\phantom{br}$

In particular, there must exist a functor 
$K(\mathcal{C}, R):\mathcal{X}\longrightarrow\mathcal{C}^{RL}$ in order to define this component.\\

On objects, this functor is defined by  
$K(\mathcal{C}, R)X:= (RX, R\zeta^{RX}(\sim)) = (RX,(\sim)^{RX})$. 
On morphisms, $K(\mathcal{C}, R)(w) = Rw$.\\

The object $(RX,(\sim)^{RX})$ is a $T$-algebra due to the general proof given in Proposition
\eqref{2311080057}, in particular by equations \eqref{2312110002} and \eqref{2312110102}.\\

The morphism $Rw:RX\longrightarrow RY$ is a morphism of $T$-algebras by taking into account the 
procedure \eqref{2312122247} with the right-hand side $Rm\cdot v$, with  $v:A\longrightarrow RX$, from 
which it can be proved in particular that $(Rm\cdot v)^{RY} = Rm\cdot v^{RX}$.\\

Finally, the definition for the functor $K(\mathcal{C}, R)$ make the associated right functor diagram
commutative \eqref{2312222042} hence the associated natural transformation is the identity 
$\rho^{KR} = 1_{R}$.\\

\subsection{The Counit}\label{2312172347}

The counit of the 2-adjunction has to have the following type

\begin{equation}\label{2311230201}
\varepsilon^{\Phi\Psi}: \Phi\Psi\longrightarrow 1_{\mathbf{EMnd}}
:\mathbf{EMnd}(Cat)\longrightarrow\mathbf{EMnd}(Cat)
\end{equation}

$\phantom{br}$

In order to propose the components, consider an extensive monad 
$(\mathcal{C}, S, \eta^{S}, (\sim)^{S})$ then its image through $\Phi\Psi$ is 
$\Phi\Psi(\mathcal{C}, S, \eta^{S}, (\sim)^{S}) = \Phi(\mathcal{C}, U^{S}, \eta^{UFS}, \zeta^{US}(\sim))
= (\mathcal{C}, U^{S}F^{S}, \eta^{UFS}, U^{S}\zeta^{UFS}(\sim)) = (\mathcal{C}, S, \eta^{S}, (\sim)^{S})$.\\

The last equality is justified by the following sequence of statements. First, 
$U^{S}F^{S}A = U^{S}(SA,(\sim)^{SA}) = SA$. Second, by the triangle given by \eqref{2312180021}, 
$\eta^{UFS} = \eta^{S}$. It only remains to prove that $U^{S}\zeta^{UFS}(\sim) = (\sim)^{S}$.\\
On definition \eqref{2312152122}, $\zeta^{USN}(v) = v^{N}$ for 
$v:A\longrightarrow U^{S}(N, (\sim)^{N})$
then if $(N, (\sim)^{N}) = F^{S}B$, $\zeta^{UFSB}(v) = v^{SB}$. 
Therefore, $U^{S}\zeta^{UFS}(\sim) = (\sim)^{S}$ as required.\\

The component of the counit can be defined as the identity.

\begin{equation}\label{2312180041}
\varepsilon^{\Phi\Psi}(\mathcal{C}, S, \eta^{S}, (\sim)^{S}) 
:= (1_{\mathcal{C}}, (\sim)^{S}):(\mathcal{C}, S, \eta^{S}, (\sim)^{S})
\longrightarrow(\mathcal{C}, S, \eta^{S}, (\sim)^{S})
\end{equation}

$\phantom{br}$

The fulfillment of the triangular identities is given as follows.\\

\newtheorem{2312042003}{Proposition}[section]
\begin{2312042003}\label{2312042003}
The triangular identities hold for the unit $\eta^{\Psi\Phi}$ and the 
 counit $\varepsilon^{\Phi\Psi}$ defined as above.
\end{2312042003}

\begin{flushleft}
\emph{Proof}:
\end{flushleft}

First, the left triangular identity given by 
$\varepsilon^{\Phi\Psi}\Phi\circ\Phi\eta^{\Psi\Phi} = 1_{\Phi}$ is proved. 
In order to accomplish this, consider a universal arrow 
$(\mathcal{C}, R, \eta^{RL}, \zeta^{R}(\sim))$. 
If the involved parts in the identity are calculated separately, then 
$\Phi\eta^{\Psi\Phi}(\mathcal{C}, R, \eta^{RL}, \zeta^{R}(\sim)) 
= \Phi(1_{\mathcal{C}},K(R,\eta^{RL}),\rho^{K}) = (1_{\mathcal{C}}, 1_{RL})$.\\

Similarly, 
$\varepsilon^{\Phi\Psi}\Phi(\mathcal{C}, R, \eta^{RL}, \zeta^{R}(\sim)) 
= \varepsilon^{\Phi\Psi}(\mathcal{C}, RL, \eta^{RL}, (\sim)^{RL}) = (1_{\mathcal{C}}, 1_{RL})$, therefore
the composition is the identity and 
$1_{\Phi}(\mathcal{C}, R, \eta^{RL}, \zeta^{R}(\sim))$ 
$= 1_{\Phi(\mathcal{C}, R, \eta^{RL}, \zeta^{R}(\sim))} = 1_{(\mathcal{C}, RL)} = (1_{\mathcal{C}}, 1_{RL})$.\\

Second, the right triangular identity given by 
$\Psi\varepsilon^{\Phi\Psi}\circ\eta^{\Psi\Phi}\Psi = 1_{\Psi}$ is proved. Consider an extensive monad
$(\mathcal{C},S,\eta^{S},(\sim)^{S})$. Let us calculate the identity by parts, for example, 
$\eta^{\Psi\Phi}\Psi(\mathcal{C},S,\eta^{S},(\sim)^{S})=$
$\eta^{\Psi\Phi}(\mathcal{C},U^{S},\eta^{U\!F\!S},$ $\zeta^{U\!S}(\sim))$. Before continue any further, 
let us remember that

\begin{equation*}
\Phi(\mathcal{C},U^{S},\eta^{U\!F\!S}, \zeta^{U\!S}(\sim)) 
= (\mathcal{C},U^{S}\!F^{S},\eta^{U\!F\!S}, U^{S}\zeta^{U\!F\!S}(\sim)) 
= (\mathcal{C},S,\eta^{S},(\sim)^{S})
\end{equation*}

$\phantom{br}$

Then 
$\eta^{\Psi\Phi}(\mathcal{C},U^{S},\eta^{U\!F\!S}, \zeta^{U\!S}(\sim))$ 
$=$ $(1_{\mathcal{C}}, 1_{\mathcal{C}^{S}},\rho^{KU\!S})$.\\

On the other hand, 
$\Psi\varepsilon^{\Phi\Psi}(\mathcal{C},S,\eta^{S},(\sim)^{S})$ $=$ 
$\Psi(1_{\mathcal{C}}, 1_{S}) = (1_{\mathcal{C}}, 1_{\mathcal{C}^{S}},\rho^{KU\!S})$. 
This concludes the identity since 
$1_{\Psi}(\mathcal{C},S,\eta^{S},(\sim)^{S}) = (1_{\mathcal{C}}, 1_{\mathcal{C}^{S}},\rho^{KU\!S})$.

\begin{flushright}
$\Box$
\end{flushright}

Consider a 2-natural transformation 
$\Gamma:F\longrightarrow G:\mathcal{A}\longrightarrow\mathcal{N}$. One of the conditions for 
$\Gamma$ to be a 2-natural transformation is to fulfill the following condition for a 2-cell
in $\mathcal{A}$. If $\alpha:f\longrightarrow g:A\longrightarrow B$ then 
$G\alpha\cdot\Gamma A = \Gamma B\cdot F\alpha$. This condition is translated for the unit,  
$\eta^{\Psi\Phi}:$
$1_{\mathbf{UArr}(Cat)}\longrightarrow \Psi\Phi:\mathbf{UArr}(Cat)\longrightarrow\mathbf{UArr}(Cat)$
for 
$(\alpha,\beta):$
$(J,V,\rho^{JV})\longrightarrow(K,W,\rho^{KW}):(\mathcal{C},R)\longrightarrow(\mathcal{D},R^{\dagger})$, 
into

\begin{equation*}
\widehat{\alpha}\cdot K(L\dashv R) = K(L^{\dagger}\dashv R^{\dagger})\cdot\widehat{\beta}
\end{equation*}

Furthermore, the condition imposed for $X$, in $\mathcal{X}$, is $\alpha RX = R^{\dagger}\beta X$ which is 
not the case  but \eqref{2311080011}. This tells us that the unit fails to be a 2-natural transformation. 
Nonetheless, as soon as one impose the commutativity condition on 1-cells for 
\eqref{2312222042}, there exists a proper 2-adjunction $\Phi\dashv \Psi$.\\

\textbf{Note}: The commutativity condition does not affect any of the proofs. That is to say,
they all can be seen as mere corollaries of the construction given so far.

\newtheorem{2312310127}[2312042003]{Theorem}
\begin{2312310127}\label{2312310127}
There is a 2-adjunction of type 
$\Phi\dashv \Psi:\mathbf{UArr}(Cat)\longrightarrow\mathbf{EMnd}(Cat)$.
\end{2312310127}
 
\section{Isomorphisms between 2-categories}\label{2312180044}

In this section, a pair of 2-isomorphisms is constructed. The first one relates, on 
$0$-cells, universal arrows and adjunctions. That is to say, 
$\mathbf{UArr}(Cat) \cong \mathbf{Adj}_{R}(Cat)$. 
The second relates, on $0$-cells, extensive monads and classical monads. That is to say, 
$\mathbf{EMnd}(Cat) \cong \mathbf{Mnd}_{E}(Cat)$.\\

\subsection{The isomorphism between $\mathbf{UArr}(Cat)$ and $\mathbf{Adj}_{R}(Cat)$}\label{2311131224}

Let us construct the 2-functor $F:\mathbf{UArr}(Cat)\longrightarrow\mathbf{Adj}_{R}(Cat)$. 
On $0$-cells, given a universal arrow $(\mathcal{C}, R,\eta^{RL},$ $\zeta^{R}(\sim))$ 
an adjunction $L\dashv R:\mathcal{C}\longrightarrow\mathcal{X}$ has to be constructed.\\

First, the function on objects 
$L:\textrm{Obj}(\mathcal{C})\longrightarrow\textrm{Obj}(\mathcal{X})$ has to
be extended to a functor. Therefore, for 
$f:A\longrightarrow B$, $Lf:= \zeta^{RLB}(\eta^{RL}B\cdot f)$. Second,
the unit of the universal arrow is proposed to be the unit of the adjunction.\\

Finally, the counit can be defined as

\begin{equation}\label{2312182230}
\varepsilon^{LR}X := \zeta^{RX}(1_{RX})
\end{equation}

\noindent for $X$ in $\mathcal{X}$.\\

For example, the triangular identity $R\varepsilon^{LR}X\cdot \eta^{RL}\!RX = 1_{RX}$, 
on $X$ in $\mathcal{X}$, is obtained from the commutativity of the following universal 
triangle.\\

\begin{equation*}
\xy<1cm,1cm>
\POS (16,27) *+{RX} = "a12",
\POS (0, 0) *+{RLRX} = "a21",
\POS (32, 0) *+{RX} = "a23", 
\POS "a12" \ar_{\eta^{RL}\!RX} "a21",
\POS "a12" \ar^{1_{RX}} "a23",
\POS "a21" \ar^{R\zeta^{RX}(1_{RX})}_{R\varepsilon^{LR}X} "a23"
\endxy
\end{equation*}

$\phantom{br}$

Therefore, 

\begin{equation*}
F(\mathcal{C}, R,\eta^{RL},\zeta^{R}(\sim)) 
= (L\dashv R:\mathcal{C}\longrightarrow\mathcal{X}, \eta^{RL}, \varepsilon^{LR})
\end{equation*}

$\phantom{br}$

For morphisms, consider 
$(J, V, \rho^{JV}):(\mathcal{C}, R, \eta^{RL}, \zeta^{R}(\sim))\longrightarrow 
(\mathcal{D}, R^{\dagger}, \eta^{RL\dagger}, \zeta^{R\dagger}(\sim))$.\\

In order to get a morphism of adjunctions, cf \cite{loj-apk}, a natural transformation of the
following type has to be defined 
$\lambda^{JV}:L^{\dagger}J\longrightarrow VL:\mathcal{C}\longrightarrow\mathcal{Y}$ 
and show that it is the \emph{adjoint dual} of $\rho^{JV}$.

\newtheorem{2311260054}{Proposition}[subsection]
\begin{2311260054}\label{2311260054}
There exists a natural transformation of type 
$\lambda^{JV}:L^{\dagger}J\longrightarrow VL:\mathcal{C}\longrightarrow\mathcal{Y}$ and its 
\emph{dual adjoint} is $\rho^{JV}$.
\end{2311260054}

\begin{flushleft}
\emph{Proof}:
\end{flushleft}

Consider the following universal triangle

\begin{equation}\label{2311262034}
\xy
\POS (25,32) *+{JA} = "a12",
\POS (37,16) *+{JRLA} = "a23",
\POS (0, 0) *+{R^{\dagger}L^{\dagger}JA} = "a31",
\POS (50, 0) *+{R^{\dagger}VLA} = "a34", 
\POS "a12" \ar_{\eta^{RL^{\dagger}\!JA}} "a31",
\POS "a31" \ar_-{R^{\dagger}\zeta^{R\dagger V\!LA}(\rho^{JV}\!LA\cdot J\eta^{RL}\!A)} "a34"
\POS "a12" \ar^{ J\eta^{RL}\!A} "a23",
\POS "a23" \ar^{\rho^{JV}\!LA} "a34",
\endxy
\end{equation}

$\phantom{br}$

Then define $\lambda^{JV}A := \zeta^{R\dagger V\!LA}(\rho^{JV}LA\cdot J\eta^{RL}\!A)$. 
Now consider $f:A\longrightarrow B$ in $\mathcal{C}$.

\begin{equation*}
\begin{array}{ccc}
\xy
\POS (25,32) *+{JA} = "a12",
\POS (54,32) *+{JB} = "a15",
\POS (37,16) *+{JRLA} = "a23",
\POS (66,16) *+{JRLB} = "a26",
\POS (0, 0) *+{R^{\dagger}L^{\dagger}JA} = "a31",
\POS (50, 0) *+{R^{\dagger}VLA} = "a34",
\POS (79, 0) *+{R^{\dagger}VLB} = "a37",
\POS (37,-17) *+{R^{\dagger}L^{\dagger}JB} = "a43",
\POS (74,-17) *+{\scriptstyle{R^{\dagger}\zeta^{R\dagger V\!LB}(\rho^{JV}\!\!LB\cdot J\eta^{RL}\!B)}} = "z",
\POS "a12" \ar_{\eta^{RL^{\dagger}}\!JA} "a31",
\POS "a31" \ar@/_1.4pc/_{R^{\dagger}L^{\dagger}\!JB} "a43",
\POS "a43" \ar@/_1.4pc/_{} "a37",
\POS "a31" \ar_-{R^{\dagger}\zeta^{R\dagger V\!LA}(\rho^{JV}\!LA\cdot J\eta^{RL}\!A)} "a34"
\POS "a12" \ar^{ J\eta^{RL}\!A} "a23",
\POS "a23" \ar^{\rho^{JV}\!LA} "a34",
\POS "a34" \ar_{R^{\dagger}VLf} "a37",
\POS "a23" \ar_{JRLf} "a26",
\POS "a26" \ar^{\rho^{JV}\!LB} "a37",
\POS "a12" \ar^{Jf} "a15",
\POS "a15" \ar^{ J\eta^{RL}\!B} "a26",
\endxy & 
\xy
\POS (0,32) *+{\phantom{0}} = "a11",
\POS (0,0) *+{\phantom{0}} = "a21",
\POS (0,-1) *+{,} = "a31",
\endxy
&
\xy
\POS (25,32) *+{JB} = "a12",
\POS (37,16) *+{JRLB} = "a23",
\POS (0, 0) *+{R^{\dagger}L^{\dagger}JB} = "a31",
\POS (50, 0) *+{R^{\dagger}VLB} = "a34", 
\POS "a12" \ar_{\eta^{RL^{\dagger}\!JB}} "a31",
\POS "a31" \ar_-{R^{\dagger}\zeta^{R\dagger V\!LB}(\rho^{JV}\!LB\cdot J\eta^{RL}\!B)} "a34"
\POS "a12" \ar^{ J\eta^{RL}\!B} "a23",
\POS "a23" \ar^{\rho^{JV}\!LB} "a34",
\endxy
\end{array}
\end{equation*}

$\phantom{br}$

The previous left diagram resembles the procedure used before. In this diagram one can 
compare a small universal triangle with right-hand side given by $\rho^{JV}LA\cdot J\eta^{RL}A$, 
which gives the component $\lambda^{JV}\!A$, with the big universal triangle with right-hand 
side $\rho^{JV}\!LB\cdot J\eta^{RL}\!B\cdot Jf$.\\

The two squares involved in the referred comparison of universal triangles commute because 
of the naturality of $J\eta^{RL}$ over the morphism $f$ and the naturality of $\rho^{JV}\!L$ 
over the same morphism.\\

Note also that the right diagram can be inserted in the left big universal diagram with 
right-hand side $\rho^{JV}\!LB\cdot J\eta^{RL}\!B\cdot Jf$. Out of this comparisons one obtains 
two parallel $R^{\dagger}$-bases making commutative this universal triangle. Therefore, this 
diagram commutes on its own and gives the naturality of $\lambda^{JV}$ over the morphism $f$ as 
required.\\

Finally, in order to prove that $\lambda^{JV}$ is the dual adjoint of $\rho^{JV}$, the following 
equation has to be fulfilled 
$\rho^{JV} = R^{\dagger}V\varepsilon^{LR}\circ R^{\dagger}\lambda^{JV}\!R\circ \eta^{RL\dagger}JR$. Since 
it can be proved that the triangular identities hold, the previous equation is equivalent to the 
following one $J\eta^{RL}\circ \rho^{JV}\!L = R^{\dagger}\lambda^{JV}\!\circ \eta^{RL\dagger}J$, but 
this is the commutativity property of the universal triangle in \eqref{2311262034}.\\ 

\begin{flushright}
$\Box$
\end{flushright}

Therefore, for morphisms

\begin{equation*}
F(J, V, \rho^{JV}) = (J, V, \rho^{JV}, \lambda^{JV})
\end{equation*}

\noindent where $\lambda^{JV}$ is the dual adjoint for $\rho^{JV}$.\\

Now consider, 2-cells in the 2-category $\mathbf{UArr}(Cat)$

\begin{equation*}
(\alpha,\beta):
(J, V, \rho^{JV})\longrightarrow (K, W, \rho^{KW}):
(\mathcal{C}, R, \eta^{RL}, \zeta^{R}(\sim))\longrightarrow
(\mathcal{D}, R^{\dagger}, \eta^{RL\dagger}, \zeta^{R\dagger}(\sim))
\end{equation*}

If this pair were to be extended to a morphism of adjunctions then there remains 
to prove the following proposition.

\newtheorem{2311262212}[2311260054]{Proposition}
\begin{2311262212}\label{2311262212}
The following equation can be derived

\begin{equation*}
\lambda^{KW}\!\circ L^{\dagger}\alpha = \beta L\circ\lambda^{JV}
\end{equation*}

\end{2311262212}

\noindent \emph{Proof}:

\begin{equation*}
\xy<1cm,1cm>
\POS (25,32) *+{JA} = "a12",
\POS (54,32) *+{KA} = "a15",
\POS (37,16) *+{JRLA} = "a23",
\POS (66,16) *+{KRLA} = "a26",
\POS (0, 0) *+{R^{\dagger}L^{\dagger}JA} = "a31",
\POS (50, 0) *+{R^{\dagger}VLA} = "a34",
\POS (79, 0) *+{R^{\dagger}WLA} = "a37",
\POS (37,-17) *+{R^{\dagger}L^{\dagger}KA} = "a43",
\POS (77,-16) *+{\scriptstyle{R^{\dagger}\zeta^{R\dagger W\!LA}(\rho^{KW}\!\!LA\cdot K\eta^{RL}\!A)}} = "z",
\POS "a12" \ar_{\eta^{RL^{\dagger}\!JA}} "a31",
\POS "a31" \ar@/_1.4pc/_{R^{\dagger}L^{\dagger}\!\alpha A} "a43",
\POS "a43" \ar@/_1.4pc/_{} "a37",
\POS "a31" \ar_-{R^{\dagger}\zeta^{R\dagger V\!LA}(\rho^{JV}\!LA\cdot J\eta^{RL}\!A)} "a34"
\POS "a12" \ar^{ J\eta^{RL}\!A} "a23",
\POS "a23" \ar^{\rho^{JV}\!LA} "a34",
\POS "a34" \ar_{R^{\dagger}\beta LA} "a37",
\POS "a23" \ar_{\alpha RLA} "a26",
\POS "a26" \ar^{\rho^{KW}\!LA} "a37",
\POS "a12" \ar^{\alpha A} "a15",
\POS "a15" \ar^{K\eta^{RL}\!A} "a26",
\endxy
\end{equation*}

$\phantom{br}$

The argument is the same as in the previous Proposition.

\begin{flushright}
$\Box$
\end{flushright}

Therefore, $F(\alpha,\beta) = (\alpha, \beta)$.\\

Let us construct the inverse 2-functor to $F$. The definition of this 2-functor is easier since
it forgets structure rather than constructing it. The 2-functor 
$G:\mathbf{Adj}_{R}(Cat)\longrightarrow\mathbf{UArr}(Cat)$, acts as follows. 
On objects, $G(L\dashv R, \eta^{RL}\!, \varepsilon^{LR}) := (\mathcal{C}, R, \eta^{RL}\!, \zeta^{R}(\sim))$ 
where the induced universal function is defined as follows

\begin{equation}\label{2311262246}
\zeta^{RX}(v) := \varepsilon^{LR}X\cdot Lv 
\end{equation}

$\phantom{br}$

On 1-cells, for $(J,V,\rho^{JV},\lambda^{JV}):L\dashv R\longrightarrow L^{\dagger}\dashv R^{\dagger}$, 
$G(J,V,\rho^{JV},\lambda^{JV}) := (J, V,\rho^{JV})$, without making any reference to the dual adjoint.\\

On 2-cells, for 
$(\alpha,\beta):(J,V,\lambda^{JV},\rho^{JV})\longrightarrow(K,W,\lambda^{KW},\rho^{KW}):
L\dashv R\longrightarrow L^{\dagger}\dashv R^{\dagger}$, $G(\alpha,\beta) := (\alpha, \beta)$.

$\phantom{br}$

\newtheorem{2311262343}[2311260054]{Proposition}
\begin{2311262343}\label{2311262343}
The 2-functors $F$ and $G$ are each other inverses.
\end{2311262343}

\begin{flushleft}
\emph{Proof}:
\end{flushleft}

One non-trivial issue about this proof is to show that $GF(\zeta^{R}(\sim)) = \zeta^{R}(\sim)$. 
In order to do it, consider
$GF(\zeta^{R}(\sim))^{(X)}(v) = G(\zeta^{RX}(1_{RX}))(v) = \zeta^{RX}(1_{RX})\cdot Lv$
using \eqref{2312182230} and \eqref{2311262246}. \\

The equality $\zeta^{RX}(1_{RX})\cdot Lv = \zeta^{RX}(v)$ has to be proved.
As before, if one draws a universal triangle
with right-hand side $\eta^{RL}RX\!\cdot v$, then
$R\zeta^{RLRX}(\eta^{RL}RX\cdot v)\cdot\eta^{RL}\!A = \eta^{RL}RX\cdot v$
and if this is pulled forward by $R\zeta^{RX}(1_{RX})$

\begin{eqnarray*}
R\zeta^{RX}(1_{RX})\cdot R\zeta^{RLRX}(\eta^{RL}RX\cdot v)\cdot\eta^{RL}\!A  
&=& R\zeta^{RX}(1_{RX})\cdot \eta^{RL}RX\cdot v\\
&=& R\zeta^{RX}(1_{RX})\cdot \eta^{RL}RX\cdot R\zeta^{RX}(v)\cdot\eta^{RL}\!A\\
&=& R\zeta^{RX}(v)\cdot\eta^{RL}\!A
\end{eqnarray*}

In the second equality, the universal triangle with right-hand side $v$ was used.
In the third equality, one of the triangular identities was applied. Now, since both 
$R$-bases make the same universal triangle commutative then 
$\zeta^{RX}(1_{RX})\cdot Lv = \zeta^{RX}(v)$ and $GF(\zeta^{R}(\sim)) = \zeta^{R}(\sim)$
as required.\\ 

The second non-trivial issue is the following one 
$FG(\varepsilon)(X) = F(\varepsilon^{LR}X\cdot L(\sim)) = \varepsilon^{LR}X\cdot L(1_{RX}) 
= \varepsilon^{LR}X$.\\

A third and final non-trivial fact is that $FG(\lambda^{JV}) = \lambda^{JV}$ but since dual adjoints are 
unique this has to hold.  

\begin{flushright}
$\Box$
\end{flushright}

\subsection{The isomorphism between $\mathbf{Mnd}_{E}(Cat)$ and $\mathbf{EMnd}(Cat)$}\label{2311131223}

In this subsection the 2-functor $H:\mathbf{EMnd}(Cat)\longrightarrow\mathbf{Mnd}_{E}(Cat)$ 
is constructed. \\

On 0-cells, that is to say extensive monads, a classical monad has to be 
given. Consider an extensive monad $(\mathcal{C}, S, \eta^{S}, (\sim)^{S})$. First, the function 
$S:\textrm{Obj}(\mathcal{C})\longrightarrow \textrm{Obj}(\mathcal{C})$ has
to be extended to a endo-functor, in order to do it, take a morphism of type $f:A\longrightarrow B$
then $Sf:= (\eta^{S}B\cdot f)^{SB}$.\\

The naturality of the universal arrow $\eta^{S}$ is proved by the following equation
$Sf\cdot\eta^{S}A = (\eta^{S}B\cdot f)^{SB}\cdot\eta^{S}A = \eta^{S}B\cdot f$.
Finally, $\mu^{S} A: = (1_{SA})^{SA}$ which will result in the multiplication of the
monad.\\

For 1-cells or morphisms, consider 
$(P,(\sim)^{PS}):(\mathcal{C}, S, (\sim)^{S})\longrightarrow(\mathcal{D}, T, (\sim)^{T})$, then
$H(P,(\sim)^{PS}) = (P,\varphi^{P})$. The natural transformation 
$\varphi^{P}:TP\longrightarrow PS:\mathcal{C}\longrightarrow\mathcal{D}$ is defined as follows

\begin{equation}\label{2311121817}
\varphi^{P}A := (P\eta^{S}\!A)^{PSA}: TPA\longrightarrow PSA
\end{equation}

$\phantom{br}$

The proof for the compatibility with the monad units is the equation 
$(P\eta^{S}\!A)^{PSA}\cdot\eta^{T}PA = P\eta^{S}\!A$. \\

For the compatibility with 
the multiplications of the monads, consider the equation 
$\varphi^{P}\circ\mu^{T}P = P\mu^{S}\circ \varphi^{P}S\circ T\varphi^{P}$  and use the definition
coming from the extensive monad. Therefore the equation to be solved is the following one
$(P\eta^{S}A)^{PSA}\cdot (1_{TPA})^{TPA} 
= P(1_{SA})^{SA}\cdot (P\eta^{S}SA)^{PSSA}\cdot [\:\eta^{T}PSA\cdot(P\eta^{S}A)^{PSA}]^{PSA}$. 

\begin{eqnarray*}
(P\eta^{S}A)^{PSA}\cdot (1_{TPA})^{TPA} 
&=& [(P\eta^{S}A)^{PSA}\cdot 1_{TPA}]^{PSA} \\
&=& [(P\eta^{S}A)^{PSA}]^{PSA}\\
&=& [P1_{SA}\cdot (P\eta^{S}A)^{PSA}]^{PSA}\\
&=& [P1_{SA}^{\:SA}\cdot P\eta^{S}SA\cdot (P\eta^{S}A)^{PSA}]^{PSA}\\
&=& P1_{SA}^{\:SA}\cdot [P\eta^{S}SA\cdot (P\eta^{S}A)^{PSA}]^{PSSA}\\
&=& P1_{SA}^{\:SA}\cdot [(P\eta^{S}SA)^{PSSA}\cdot \eta^{T}PSA\cdot(P\eta^{S}A)^{PSA}]^{PSSA}\\
&=& P1_{SA}^{\:SA}\cdot (P\eta^{S}SA)^{PSSA}\cdot [\eta^{T}PSA\cdot(P\eta^{S}A)^{PSA}]^{PSA} 
\end{eqnarray*}

$\phantom{br}$

In the first equality, the fact that $(P\eta^{S}A)^{PSA}$ is a $T$-algebraic morphism was used.
In the fourth equality, the property $1_{SA} = (1_{SA})^{SA}\cdot\eta^{S}SA$ was applied.
In the fifth equality, the fact that $P1_{SA}^{\:SA}$ is a morphism of $T$-algebras was used.
In the sixth equality, the property $P\eta^{S}SA = (P\eta^{S}SA)^{PSSA}\cdot \eta^{T}PSA$ was applied.
Finally, in the seventh equality, that $(P\eta^{S}A)^{PSA}$ is a $T$-algebraic morphism was used.\\ 

On 2-cells, a transformation of extensive monads 
$\theta:(P,(\sim)^{P})\longrightarrow(Q,(\sim)^{Q}):
(\mathcal{C}, S, (\sim)^{S})\longrightarrow(\mathcal{D}, T, (\sim)^{T})$ has to be taken into a
transformation of classical monads. Then the following property
$\psi^{Q}\circ T\theta = \theta S\circ\varphi^{P}$ has to be proved. By taking into account 
the respective definitions, this requirement is translated into 
$(Q\eta^{S}\!A)^{QSA}\cdot(\eta^{T}QA\cdot\theta A)^{TQA}= \theta SA\cdot(P\eta^{S}\!A)^{PSA}$.\\

\begin{eqnarray*}
(Q\eta^{S}\!A)^{QSA}\cdot(\eta^{T}QA\cdot\theta A)^{TQA} 
&=& \big((Q\eta^{S}\!A)^{QSA}\cdot\eta^{T}QA\cdot\theta A \big)^{QSA}\\
&=& (Q\eta^{S}\!A \cdot\theta A)^{QSA}\\
&=& (\theta SA\cdot P\eta^{S}\!A)^{QSA}\\
&=& \theta SA\cdot(P\eta^{S}\!A)^{PSA}
\end{eqnarray*}

In the first equality, the fact that $(Q\eta^{S}\!A)^{QSA}$ is an algebraic $T$-morphism is used.
In the second equality, the property $(Q\eta^{S}\!A)^{QSA}\cdot\eta^{T}QA = Q\eta^{S}\!A$ was applied.
In the third equality, the naturality of $\theta$ over the morphism $\eta^{S}\!A$ was used.
In the fourth equality, that the components for $\theta S$ are morphisms of $T$-algebras was applied.\\

This finishes the construction for the 2-functor 
$H:\mathbf{EMnd}(Cat)\longrightarrow\mathbf{Mnd}_{E}(Cat)$.\\

For the 2-functor $K:\mathbf{Mnd}_{E}(Cat)\longrightarrow\mathbf{EMnd}(Cat)$, 
consider the 0-cells and take a classic monad $(\mathcal{C},S,\eta^{S},\mu^{S})$ then it only 
remains to define the associated extensive function. Consider a morphism of type $h:A\longrightarrow SB$,
in $\mathcal{C}$, the extensive function is defined as follows.

\begin{equation}\label{2311270042}
h^{SB} := \mu^{S}B\cdot Sf
\end{equation}

$\phantom{br}$

Out of this definition, the property $k^{SC}\cdot h^{SB} = (k^{SC}\cdot h)^{SC}$ has to be fulfilled.

\begin{eqnarray*}
\mu^{S}C\cdot Sk\cdot \mu^{S}B\cdot Sh
&=& \mu^{S}C\cdot \mu^{S}SC \cdot SSk\cdot Sh\\
&=& \mu^{S}C\cdot S\mu^{S}C \cdot SSk\cdot Sh\\
&=& \mu^{S}C\cdot S(\mu^{S}C\cdot Sk\cdot h)\\
\end{eqnarray*}

$\phantom{br}$

In the first equality, the naturality of $\mu^{S}$ over the morphism $k$ was used.
In the second equality, the associativity of the multiplication $\mu^{S}$ was applied. Then

\begin{equation}\label{2312191920}
K(\mathcal{C},S,\eta^{S},\mu^{S}) = (\mathcal{C},S,\eta^{S},(\sim)^{S})
\end{equation}

$\phantom{br}$

On 1-cells, for a morphism of monads 
$(P,\varphi^{P}):(\mathcal{C},S)\longrightarrow(\mathcal{D},T)$, 
the function for the extensive algebra $PSA$, in $\mathcal{D}$, has to be defined.

\begin{equation}\label{2311270047}
r^{PSA} := P\mu^{S}\!A\cdot\varphi^{P}\!SA\cdot Tr
\end{equation}

$\phantom{br}$

The property for $(PSA, (\sim)^{PSA})$ to be a $T$-algebra is the following one 
$r^{PSA}\!\cdot d^{TD} = (r^{PSA}\!\cdot d)^{PSA}$, where $r:D\longrightarrow PSA$.\\

\begin{eqnarray*}
r^{PSA}\!\cdot d^{TD} 
&=& P\mu^{S}\!A\cdot\varphi^{P}\!SA\cdot Tr\cdot\mu^{T}\!D\cdot Td\\
&=& P\mu^{S}\!A\cdot\varphi^{P}\!SA\cdot\mu^{T}PSA\cdot TTr\cdot Td\\
&=& P\mu^{S}\!A\cdot P\mu^{S}SA\cdot\varphi^{P}\!SSA\cdot T\varphi^{P}\!SA\cdot TTr\cdot Td\\
&=& P\mu^{S}\!A\cdot PS\mu^{S}\!A\cdot\varphi^{P}\!SSA\cdot T\varphi^{P}\!SA\cdot TTr\cdot Td\\
&=& P\mu^{S}\!A\cdot \varphi^{P}\!SA\cdot TP\mu^{S}\!A\cdot T\varphi^{P}\!SA\cdot TTr\cdot Td\\
&=& P\mu^{S}\!A\cdot \varphi^{P}\!SA\cdot T(r^{PSA}\!\cdot d)\\
&=& (r^{PSA}\!\cdot d)^{PSA}
\end{eqnarray*}

The details for the chain of equalities are explained. 
In the second equality, the naturality of $\mu^{T}$ over $r$ was used.
In the third equality, the compatibility of $\varphi^{P}$ with the multiplication of the monads was used. 
In the fourth equality, the associativity of the monad was used.
In the fifth one, the naturality of $\varphi^{P}$ over the morphism $\mu^{S}\!A$ was used.
Finally, in the last two equalities the definition of the function $(\sim)^{PSA}$ was used.\\

The property for $Ph^{SB}:(PSA, (\sim)^{PSA})\longrightarrow(PSB, (\sim)^{PSB})$ to be a morphism of 
$T$-algebras has also to be proved, for any $h:A\longrightarrow SB$, in $\mathcal{C}$.
In order to prove this statement, for $r:D\longrightarrow PSA$, 
the following equation must hold $Ph^{SB}\!\cdot r^{PSA} = (Ph^{SB}\!\cdot r)^{PSB}$.\\

\begin{eqnarray*}
Ph^{SB}\cdot r^{PSA} 
&=& P\mu^{S}B\cdot PSh\cdot P\mu^{S}\!A\cdot\varphi^{P}\!SA\cdot Tr\\
&=& P\mu^{S}B\cdot P\mu^{S}SB\cdot PSSh\cdot\varphi^{P}\!SA\cdot Tr\\
&=& P\mu^{S}B\cdot PS\mu^{S}B\cdot PSSh\cdot\varphi^{P}\!SA\cdot Tr\\
&=& P\mu^{S}B\cdot \varphi^{P}\!SB\cdot TP\mu^{S}B\cdot TPSh\cdot Tr\\
&=& P\mu^{S}B\cdot \varphi^{P}\!SB\cdot T(P\mu^{S}B\cdot PSh\cdot r)\\
&=& P\mu^{S}B\cdot \varphi^{P}\!SB\cdot T(Ph^{SB}\!\cdot r)\\
&=& (Ph^{SB}\cdot r)^{PSB}
\end{eqnarray*}

$\phantom{br}$

The procedure is explained as follows. 
In the second equality, the naturality of $\mu^{S}$ over $h$ was used.
In the third equality, the associativity of the multiplication $\mu^{S}$ was used.
In the fourth equality, the naturality of $\varphi^{P}$ over $\mu^{S}B\cdot Sh$ was used.
In the sixth and seventh equalities, the definition of the functions $(\sim)^{SB}$  and $(\sim)^{PSB}$
were used, respectively.\\

Then for 1-cells

\begin{equation}\label{2312192002}
K(P, \varphi^{P}) = (P, (\sim)^{PS})
\end{equation}

$\phantom{br}$

On 2-cells, for the transformation of monads 
$\theta:(P,\varphi^{P})\longrightarrow(Q,\psi^{Q}):(\mathcal{C},S)\longrightarrow(\mathcal{D},T)$, 
the image for the 2-functor is $K(\theta) = \theta$. 
This image has to be a morphism of $T$-algebras of type 
$\theta SA:(PSA, (\sim)^{PSA})\longrightarrow(QSA, (\sim)^{QSA})$, then for $r:D\longrightarrow PSA$
it has to be proved that $\theta SA\cdot r^{PSA} = (\theta SA\cdot r)^{QSA}$.

\begin{eqnarray*}
\theta SA\cdot r^{PSA} 
&=& \theta SA\cdot P\mu^{S}\!A \cdot \varphi^{P}\!SA \cdot Tr\\
&=& Q\mu^{S}\!A\cdot \theta SSA \cdot \varphi^{P}\!SA \cdot Tr\\
&=& Q\mu^{S}\!A\cdot \psi^{Q}\!SA \cdot T\theta SA \cdot Tr\\
&=& Q\mu^{S}\!A\cdot \psi^{Q}\!SA \cdot T(\theta SA \cdot r)\\
&=& (\theta SA\cdot r)^{QSA}
\end{eqnarray*}

In the second equality, the naturality of $\theta$ over $\mu^{S}\!A$ was used.
In the third equality, the property for $\theta$ to be a transformation of monads was used.
In the fifth equality, the definition of the function $(\sim)^{QSA}$ was used.

\newtheorem{2311270104}{Proposition}[subsection]
\begin{2311270104}\label{2311270104}
The 2-functors $H$ and $K$ are each other inverses.
\end{2311270104}

\begin{flushleft}
\emph{Proof}:
\end{flushleft}

The only non-trivial elements to be proved for the composition $KH$ is the following one,
$KH((\sim)^{PS})(A)(r) = (\sim)^{PS}(A)(r)$. First, $KH((\sim)^{PS})(A)(r) = K((P\eta^{S}A)^{PSA})(r) 
= P(1_{SA})^{SA}\cdot(P\eta^{S}A)^{PSA}\cdot(\eta^{T}PSA\cdot r)^{TPSA}$ then 

\begin{eqnarray*}
P(1_{SA})^{SA}\cdot(P\eta^{S}\!A)^{PSA}\cdot(\eta^{T}PSA\cdot r)^{TPSA} 
&=& \big(P(1_{SA})^{SA}\cdot P\eta^{S}\!A\big)^{PSA}\cdot(\eta^{T}PSA\cdot r)^{TPSA}\\
&=& (1_{PSA})^{PSA}\cdot(\eta^{T}PSA\cdot r)^{TPSA}\\
&=& ((1_{PSA})^{PSA}\cdot\eta^{T}PSA\cdot r)^{PSA}\\
&=& r^{PSA}\\
&=& (\sim)^{PS}(A)(r)
\end{eqnarray*}

In the first equality, the fact that $P(1_{SA})^{SA}$ is a morphism of $T$-algebras was applied.
In the second equality, the property $(1_{SA})^{SA}\cdot\eta^{S}\!A = 1_{SA}$ was used.
In the third equality, the fact that $(1_{PSA})^{PSA}$ is a $T$-algebraic morphism was applied. 
Finally, in the fourth equality, the property $(1_{PSA})^{PSA}\cdot\eta^{T}PSA = 1_{PSA}$ was used.\\

The only non-trivial element to be proved for the composition $HK$ is  for 1-cells. 
Then 
$HK(\varphi^{P})(A) = H(P\mu^{S}\!A\cdot\varphi^{P}\!SA\cdot T(\sim)) 
= P\mu^{S}\!A\cdot\varphi^{P}\!SA\cdot T(P\eta^{S}A)$ but 
$P\mu^{S}\!A\cdot\varphi^{P}\!SA\cdot T(P\eta^{S}A) 
= P\mu^{S}\!A\cdot PS\eta^{S}\!A\cdot\varphi^{P}\!A = \varphi^{P}\!A$ by using the naturality
of $\varphi^{P}$ and one of the unitalities of the monad.

\begin{flushright}
$\Box$
\end{flushright}

All of the calculations from this monograph can be summarized into the 
following Theorem.

\newtheorem{2312191947}[2311270104]{Theorem}
\begin{2312191947}\label{2312191947}
The following diagram commutes serially

\begin{equation}\label{2312191949}
\xy
\POS (0,20) *+{\mathbf{UArr}(Cat)} = "a11",
\POS (44,20) *+{\mathbf{Adj}_{R}(Cat)} = "a12",
\POS (0,0) *+{\mathbf{EMnd}(Cat)} = "a21",
\POS (44,0) *+{\mathbf{Mnd}(Cat)} = "a22",
\POS "a11" \ar@<-4pt>_{\Phi}^{\dashv} "a21",
\POS "a21" \ar@<-6pt>_{\Psi} "a11",
\POS "a21" \ar@<-2.5pt>_{H} "a22",
\POS "a22" \ar@<-2.5pt>_{K} "a21",
\POS "a11" \ar@<-2.5pt>_{F} "a12",
\POS "a12" \ar@<-2.5pt>_{G} "a11",
\POS "a12" \ar@<-4pt>_{\Phi_{E}}^{\dashv} "a22",
\POS "a22" \ar@<-6pt>_{\Psi_{E}} "a12",
\endxy
\end{equation}

\end{2312191947}

\section{Conclusions}

This article might serve as a 1-dimensional template in order to formulate higher pseudo adjunctions 
\cite{grj-foc} that relate higher expressions of universal arrows and extensive, or no-iteration, 
pseudo monads. For example, if one extends these calculations to a 2-categorical layout, 
2-$\mathbf{Cat}$, then one has to take any equation and transform it into an invertible $n$-cell and 
then impose coherence conditions on these higher cells.\\

The higher context for this monography has very important advances in \cite{maf-nip} for (extensive) 
pseudo monads and in \cite{fit-psl} for pseudo universal arrows. Also for the higher adjunction
relating these two concepts there are several results, \cite{gan-fot} and \cite{laa-fra}.
Actually, there remain only to check some lengthy details for such higher adjunction \cite{vaa-inp}.\\

Another important contribution of this monography is the relation between extensive monads with 
universal arrows in order to explore further properties on monadic functional programming 
\cite{moe-noc}.\\

Finally, the authors think that the relationship between universal arrows and extensive 
monads through a 2-adjunction is very important since several results only take into account 
this relationship through isomorphisms of categories of type

\begin{equation*}
\mathrm{Hom}_{\mathbf{UArr}(Cat)}((\mathcal{C},U^{S}),(\mathcal{C},U^{S}))\cong
\mathrm{Hom}_{\mathbf{EMnd}(Cat)}((\mathcal{C},S),(\mathcal{C},S))
\end{equation*}

\noindent where the objects in the first entry are called \emph{liftings}. But in this monograph
the narrative goes beyond this point since it constructs the very 2-adjunction from which all of 
these isomophisms come from.

\end{document}